\catcode`\@=11
\newif\if@fewtab\@fewtabtrue

{\count255=\time\divide\count255 by 60
\xdef\hourmin{\number\count255}
\multiply\count255 by-60\advance\count255 by\time
\xdef\hourmin{\hourmin:\ifnum\count255<10 0\fi\the\count255}}
\def\ps@draft{\let\@mkboth\@gobbletwo
    \def\@oddhead{}
    \def\@oddfoot
       {\hbox to 7 cm{$\scriptstyle Draft\ version:\ \draftdate$
       \hfil}\hskip -7cm\hfil\rm\thepage \hfil}
    \def\@evenhead{}\let\@evenfoot\@oddfoot}

\def\ceqno{\global\@fewtabfalse
    \ifcase\@eqcnt \def\@tempa{& & &}\or \def\@tempa{& &}
      \or \def\@tempa{&}
      \or\def\@tempa{}\fi\@tempa
{\rm(\theequation)}}

\def\aeqno#1{\global\@fewtabfalse
    \ifcase\@eqcnt \def\@tempa{& & &}\or \def\@tempa{& &}
      \or \def\@tempa{&}
      \or\def\@tempa{}\fi\@tempa
{\rm(\theequation,#1)}}

\def\label#1{\ifnum\draftcontrol=1
 \global\def\draftnote{$\scriptstyle #1$}\fi
 \@bsphack\if@filesw {\let\thepage\relax
   \def\protect{\noexpand\noexpand\noexpand}%
\xdef\@gtempa{\write\@auxout{\string
      \newlabel{#1}{{\@currentlabel}{\thepage}}}}}\@gtempa
   \if@nobreak \ifvmode\nobreak\fi\fi\fi
  \@esphack}

\def\alabel#1#2{\label{#1}\global\@fewtabfalse
    \ifcase\@eqcnt \def\@tempa{& & &}\or \def\@tempa{& &}
      \or \def\@tempa{&}
      \or\def\@tempa{}\fi\@tempa
{\hbox to 3cm{\phantom{\rm(\theequation,#2)}
\draftnote \hfil}\hskip -3cm {\rm(\theequation,#2)}}}

\def\clabel#1{\label{#1}\global\@fewtabfalse
    \ifcase\@eqcnt \def\@tempa{& & &}\or \def\@tempa{& &}
      \or \def\@tempa{&}
      \or\def\@tempa{}\fi\@tempa
{\hbox to 3cm{\phantom{\rm(\theequation)}
\draftnote \hfil}\hskip -3cm{\rm(\theequation)}}}

\def\eqnarray{\def\draftnote{{}}\global\@fewtabtrue
\stepcounter{equation}\let\@currentlabel=\theequation
\global\@eqnswtrue
\global\@eqcnt\z@\tabskip\@centering\let\\=\@eqncr
$$\halign to \displaywidth\bgroup\@eqnsel\hskip\@centering\@eqcnt\z@
  $\displaystyle\tabskip\z@{##}$&\global\@eqcnt\@ne
  \hskip 1\arraycolsep \hfil${##}$\hfil
  &\global\@eqcnt\tw@ \hskip 1\arraycolsep
$\displaystyle\tabskip\z@{##}$
\hfil  \tabskip\@centering&\global\@eqcnt\thr@@\llap{##}\tabskip\z@
\cr}

\def\endeqnarray{\@@eqncr\egroup
      \global\advance\c@equation\m@ne$$\global\@ignoretrue}

\def\@eqnnum{\hbox to 3cm{\phantom{\rm(\theequation)} \draftnote
                         \hfil}\hskip -3cm {\rm(\theequation)}}

\def\@@eqncr{\let\@tempa\relax
    \ifcase\@eqcnt \def\@tempa{& & &}\or \def\@tempa{& &}
      \or \def\@tempa{&}
      \or\def\@tempa{}
\fi\@tempa
\if@eqnsw
\if@fewtab\@eqnnum\fi
\stepcounter{equation}\fi\global
\@eqnswtrue\global\@eqcnt\z@\global\@fewtabtrue\cr}

\def\draftcite#1{\ifnum\draftcontrol=1#1\else{}\fi}

\def\@lbibitem[#1]#2{\item{}\hskip -3cm \hbox to 2cm
{\hfil$\scriptstyle\draftcite{#2}$}\hskip
1cm[\@biblabel{#1}]\if@filesw
     {\def\protect##1{\string ##1\space}\immediate
      \write\@auxout{\string\bibcite{#2}{#1}}}\fi\ignorespaces}

\def\@bibitem#1{\item\hskip -3cm \hbox to 2cm
{\hfil $\scriptstyle\draftcite{#1}$}\hskip 1cm
\if@filesw \immediate\write\@auxout
       {\string\bibcite{#1}{\the\value{\@listctr}}}\fi\ignorespaces}

 \def\nsection#1{\section{#1}\setcounter{equation}{0}}

     \def\nappendix#1{\vskip 1cm\no{\bf Appendix
         #1}\def\thesection{#1} \setcounter{equation}{0}}

\font\tendl=msbm10  scaled \magstep1%double line
\font\sevendl=msbm7 scaled \magstep1
\font\fivedl=msbm5 scaled \magstep1
\font\tengl=eufm10  scaled \magstep1% gothic letters
\font\sevengl=eufm7 scaled \magstep1
\font\fivegl=eufm5 scaled \magstep1

\newfam\dlfam \def\dl{\fam\dlfam\tendl} % \dl is double line
\textfont\dlfam=\tendl \scriptfont\dlfam=\sevendl
\scriptscriptfont\dlfam=\fivedl
\newfam\glfam  % \gl is gothic letters
\textfont\glfam=\tengl \scriptfont\glfam=\sevengl
\scriptscriptfont\glfam=\fivegl

\def\draftdate{\number\month/\number\day/\number\year\ \ \ \hourmin }

\global\def\draftcontrol{0}
\catcode`\@=12
\def\tilde{\widetilde}
\def\hat{\widehat}

\documentstyle[11pt,epsf]{article}

\def\theequation{{\thesection.\arabic{equation}}}
\newcommand{\be}{\begin{eqnarray}}
\newcommand{\en}{\end{eqnarray}\vs 0.5 cm}

\newcommand{\no}{\noindent}
\newcommand{\vs}{\vskip}

\newcommand{\NR}{{{\dl R}}}%letra doble raya en modo matematico
\newcommand{\NC}{{{\dl C}}}%letra doble raya en modo matematico
\newcommand{\NZ}{{{\dl Z}}}%letra doble raya en modo matematico
\newcommand{\Ng}{{{\bf g}}}
\newcommand{\Nt}{{{\bf t}}}
\newcommand{\sk}{{\sc k}}
\newcommand{\qq}{\begin{eqnarray}}

\newcommand{\ee}{{\rm e}}

\newcommand{\qqq}{\end{eqnarray}}

\newcommand{\tr}{\hbox{tr}}

\newcommand{\CA}{{\cal A}}

\newcommand{\CG}{{\cal G}}

\newcommand{\CO}{{\cal O}}

\newcommand{\CZ}{{\cal Z}}
\newcommand{\s}{\hspace{0.05cm}}
\newcommand{\m}{\hspace{0.025cm}}

\begin{document}

\title{Basic gerbe over non simply connected compact groups}

\author{\ \\Krzysztof Gaw\c{e}dzki\footnote{membre du C.N.R.S.}\\
%, Laboratoire de Physique,
%ENS-Lyon,\\46, All\'ee d'Italie, F-69364 Lyon, France\\
Nuno Reis\footnote{partially supported by the Portuguese grant Praxis XXI/BD/18138/98 from FCT and two fellowships from Funda\c{c}\~ao Calouste Gulbenkian and from French Embassy in Lisbon.}\\ 
\\Laboratoire de Physique, ENS-Lyon,\\46, All\'ee
d'Italie, F-69364 Lyon, France}
\date{ }
\maketitle

%%% for draft versions, suppress in definitive version:
%\draft
%%
%%% suppress in definite version:
%\vskip 1 cm

\begin{abstract}
\vskip 0.3cm \noindent
We present an explicit construction of the basic bundle gerbes with connection 
over all connected compact simple Lie groups. These are geometric objects that 
appear naturally in the Lagrangian approach to the WZW conformal field theories. 
Our work extends the recent construction of E. Meinrenken \cite{Meinr} 
restricted to the case of simply connected groups. 
\end{abstract}
\vskip 1.2cm

\nsection{Introduction}
\label{sec:intro}

Bundle gerbes \cite{Murray,MurrS,Chatt,Hitchin} are geometric objects 
glued from local inputs with the use of transition data forming a 1-cocycle 
of line bundles. In a version equipped with connection, they found 
application in the Lagrangian approach to string theory where they permit 
to treat in an intrinsically geometric way the Kalb-Ramond 2-form fields 
$\,B\,$ that do not exist globally \cite{Alva,Gaw0,Sharpe,GR}. 
One of the simplest situations of that type involves group manifolds 
$\,G\,$ when the (local) $\,B\,$ 
field satisfies $\,dB=H\,$ with $\,H={\sk\over 12\pi}\,\tr\,(g^{-1}dg)^3$.
Such $\,B\,$ fields appear in the WZW conformal field theories of level
$\,\sk\,$ and the related coset models \cite{WZW,GawKup}. Construction of the 
corresponding gerbes allows a systematic Lagrangian treatment of such 
models, in particular, of the conformal boundary conditions corresponding 
to open string branes. This was discussed in a detailed way in 
ref.\,\,\cite{GR}. The abstract framework was illustrated there by 
the example of the $\,SU(N)\,$ group and of groups covered by $\,SU(N)$. 
\,Here we extend the recent construction \cite{Meinr} of the basic gerbe 
on all simple, connected and simply connected compact groups to the 
non-simply connected case. Similarly as in the case of groups covered by 
$\,SU(N)$, \,this requires solving a simple cohomological equation. 
The solution exists only if the level $\,\sk\,$ is such that the closed 
3-form $\,H\,$ has periods in $\,2\pi\NZ$. \,We recover this way the 
constraints on the level first worked out in ref.\,\,\cite{FGK}, see also
\cite{KSchell}. The present construction opens the possibility to extend 
to other non-simply connected groups the classification of branes for 
the groups covered by $\,SU(N)\,$ worked out in ref.\,\,\cite{GR}.
\vskip 0.3cm

\nsection{Basic gerbe on simply connected compact groups \cite{Meinr}}

We refer the reader to ref.\,\,\cite{Murray} for an introduction to bundle 
gerbes with connection, to \cite{MurrS} for the notion of stable
isomorphisms of gerbes (employed below in accessory manner) and to 
\cite{GR} for a discussion of the relevance of the notions for the WZW 
models of conformal quantum field theory. For completeness, we shall only 
recall here the basic definition \cite{Murray}. 
\,For $\,\pi:Y\mapsto M$, \, let 
\qq
Y^{{[n]}}\ =\ Y\times_{_M}\hspace{-0.1cm}Y\dots\times_{_M}\hspace{-0.1cm}
Y\ =\ \{\,(y_1,\dots,y_n)\in Y^{n}\ |\ \pi(y_1)=\,\dots\,=\pi(y_n)\,\}
\qqq
denote the $n$-fold fiber product of $\,Y$, \,$\pi^{{[n]}}$ the obvious 
map from $Y^{{[n]}}$ to $M$ and 
$\,p_{n_1\dots n_k}$ the projection of $(y_1,\dots,y_n)$ to 
$(y_{n_1}, \dots,y_{n_k})$. \,Let $\,H\,$ be a closed 3-form on manifold
$\,M$.
\vskip 0.4cm

\noindent\hbox to 2.4cm{\bf Definition.\hfill}\parbox[t]{12.1cm}{A bundle 
gerbe $\,\CG\,$ with connection (shortly, a gerbe) of curvature $\,H\,$ over $\,M\,$ 
is a quadruple $\,(Y,B,L,\mu)$, \,where}
\vskip 0.4cm
\noindent\hbox to 3cm{\hfill 1.\ \ }\parbox[t]{11.4cm}{$Y\,$ is a manifold 
provided with a surjective
submersion $\,\pi:Y\rightarrow M$.}
\vskip 0.3cm
\noindent\hbox to 3cm{\hfill 2.\ \ }\parbox[t]{11,4cm}
{${B}\,$ is a 2-form on $Y$ such that
\qq
%\hspace{-1.1cm}
d{B}=\pi^*{H}\,,
\label{1}
\qqq}
\vskip 0.09cm
\noindent\hbox to 3cm{\hfill 3.\ \ }\parbox[t]{11.4cm}{$L\,$ is a hermitian 
line bundle with a
connection $\nabla$ over $\,Y^{{[2]}}\,$ with
the curvature form
\qq
curv(\nabla)=p_2^*{B}-p_1^*{B}\,.
\label{2}
\qqq}
\vskip 0.09cm
\noindent\hbox to 3cm{\hfill 4.\ \ }\parbox[t]{11.4cm}{$\mu:p_{12}^*{L}
\otimes p_{23}^*{L}\longrightarrow
p_{13}^*{L}\,$ is an isomorphism between the line bundles with connection
over $\,Y^{[3]}\,$ such that over $\,Y^{[4]}$ 
\qq 
\mu\circ(\mu\otimes id)\,=\,\mu\circ(id\otimes\mu)\,.
\label{assoc} 
\qqq}
\vskip 0.3cm

\noindent The 2-form ${B}$ is called the curving of the gerbe. The
isomorphism $\,\mu\,$ defines a structure of a groupoid on $\,{L}\,$ 
with the bilinear product $\,\mu:{L}_{(y_1,y_2)}
\otimes{L}_{(y_2,y_3)}\rightarrow{L}_{(y_1,y_3)}$. 
%Bundle $\,{L}\,$ restricted to the diagonal composed of the elements
%$\,(y,y)\,$ may be naturally trivialized by the choice of the
%units of the groupoid multiplication and $\,\mu\,$ determines a natural
%isomorphism between $\kappa^*{L}$ and ${L}^{-1}$, where
%$\kappa(y_1,y_2)=(y_2,y_1)$. 
\vskip 0.3cm

In ref.\,\,\cite{Meinr}, an explicit and elegant construction of a gerbe 
with curvature $\,H={1\over{12\pi}}\,\tr\,(g^{-1}dg)^3\,$ on compact, 
simple, connected and simply connected groups $\,G\,$ was given.
Since it corresponds to the lowest admissible positive level $\,\sk=1$,
such a gerbe has been termed ``basic''. It is unique up to stable
isomorphisms. We re-describe here the construction of \cite{Meinr} 
in a somewhat more concrete and less elegant terms (ref.\,\,\cite{Meinr} constructed 
the gerbe equivariant w.r.t. adjoint action; we skip here the higher order equivariant 
corrections). 
\vskip 0.3cm

Let us first collect some simple facts used in the sequel. 
\,Let $\,\Ng\,$ be the Lie algebra of $\,G\,$ and
\qq
\Ng^{\NC}\ =\ \Nt^\NC\oplus\Big(\mathop{\oplus}
\limits_{\alpha\in\Delta}\NC e_\alpha\Big)
\qqq
the root decomposition of its complexification, with $\,\Nt\,$ standing
for the Cartan algebra and $\,\Delta\,$ for the set of roots. 
Let $\,r\,$ be the rank of $\,\Ng\,$ and 
$\,\alpha_i,\ \alpha_i^{^{\hspace{-0.03cm}\vee}},\ \lambda_i,
\ \lambda^{^{\hspace{-0.03cm}\vee}}_i,\ i=1,\dots,r,$ be the simple roots,
coroots, weights and coweights of $\,\Ng$\, generating the lattices
$\,Q$, $\,Q^{^{\hspace{-0.03cm}\vee}}$, $\,P$ and 
$\,P^{^{\hspace{-0.03cm}\vee}}$, respectively. \,We shall identify
in the standard way $\,\Ng\,$ and its dual using the $ad$-invariant
bilinear form $\,\tr\,XY\,$ on $\,\Ng\,$ normalized so that the long
roots have length squared $2$. The roots and coroots satisfy 
$\,\alpha=2\alpha^{^{\hspace{-0.03cm}\vee}}/\tr
(\alpha^{^{\hspace{-0.03cm}\vee}})^2$.  \,The highest
root $\,\phi=\sum\limits_{i=1}^r k_i\alpha_i=\phi^{^{\hspace{-0.03cm}\vee}}=
\sum\limits_{i=1}^r k_i^{^{\hspace{-0.03cm}\vee}}
\alpha_i^{^{\hspace{-0.03cm}\vee}}$. \,The dual Coxeter number 
$\,h^{^{\hspace{-0.03cm}\vee}}=\sum\limits_{i=0}^r
k_i^{^{\hspace{-0.03cm}\vee}}$, \,where $\,k_0=k_0^{^{\hspace{-0.03cm}\vee}}
=1$. \,The space of conjugacy classes in $\,G$, \,i.e.\,\,of the orbits 
of the adjoint action of $\,G\,$ on itself, may be identified with 
the Weyl alcove
\qq
\CA=\{\tau\in\Nt\,\,|\,\,\tr\,\alpha_i\tau\geq 0,\ \,i=1,\dots,r,
\,\ \tr\,\phi\tau\leq 1\,\}
\qqq
since every conjugacy class has a single element of the form 
$\,\ee^{2\pi i\tau}\,$ with $\,\tau\in\CA$. \,Set $\,\CA\,$ is a simplex with 
vertices $\,\tau_i={_1\over^{k_i}}\lambda_i^{^{\hspace{-0.03cm}\vee}}\,$ 
and $\,\tau_0=0$. \,Let
\qq
\CA_0=\{\,\tau\in\CA\,\,|\,\,\tr\,\phi\tau<1\,\}
\quad{\rm and}\quad\CA_i=\{\,\tau\in\CA\,\,|\,\,\tr\,\alpha_i\tau>0\,\}\quad{\rm for}
\ \ i\not=0 
\qqq
and let $\,\CA_I=\mathop{\cap}\limits_{i\in I}\CA_i\,$
for $\,I\subset\{0,1,..\m,r\}\equiv R$, 
\,We shall denote by $\,G_i\,$ the adjoint action stabilizer
of $\,\ee^{2\pi i\tau_i}$:
\qq
G_i\ =\ \{\,\gamma\in G\,\,|\,\,\gamma\,\ee^{2\pi i\tau_i}\,\gamma^{-1}
=\ee^{2\pi i\tau_i}\,\}
\qqq
and by $\,\Ng_i\,$ its Lie algebra. The complexification of $\,\Ng_i\,$ is
\qq
\Ng_i^\NC\ =\ \Nt^\NC\oplus\Big(\mathop{\oplus}
\limits_{\alpha\in\Delta_i}\NC e_\alpha\Big)\,,
\qqq
where $\,\Delta_i\,$ is composed of roots $\,\alpha\,$ such that
$\,\tr\,\tau_i\alpha\in\NZ$. For $\,i=0$, $\,G_0=G\,$ and $\,\Ng_0=\Ng$.
\,For $\,i\not=0$, $\,\Ng_i\,$ is a simple Lie algebra with simple roots 
$\,\alpha_j,\ j\not=i,\,$ and $\,-\phi$. \,Its simple coweights are
$\,{k_j}(\tau_j-\tau_i),\ j\not=i\,$ and $\,-\tau_i\,$ and they generate 
the coweight lattice $\,P^{^{\hspace{-0.03cm}\vee}}_i\,$ of $\,\Ng_i$. 
\vskip 0.3cm

The main complication in the construction of the basic gerbe over general 
compact simply connected groups is that the stabilizers $\,G_i\,$ are 
connected but, unlike for $\,SU(N)$, \,they are not necessarily simply 
connected. We shall denote by $\,\tilde G_i\,$ their universal covers. 
$\,G_i=\tilde G_i/\CZ_i$, \,where $\,\CZ_i\,$ is the subgroup of the 
center of $\,\tilde G_i$. $\,\CZ_i\,$ is composed of elements of the 
form $\,\ee_i^{2\pi i p}\,$ with $\,p\in P^{^{\hspace{-0.03cm}\vee}}_i
\cap Q^{^{\hspace{-0.03cm}\vee}}\,$
and $\,\ee_i\,$ standing for the exponential map from $\,\Ng_i\,$
to $\,\tilde G_i$. \,Since $\,\tau_i\,$ is also a weight of $\,\Ng_i$, 
it defines a character $\,\chi_i\,$ on the Cartan subgroup 
$\,\tilde T_i\,$ of $\,\tilde G_i$, \,and hence also on 
$\,\CZ_i$, \,by the formula
\qq
\chi_i(\ee_i^{2\pi i\tau})\ =\ \ee^{2\pi i\,\tr\,\tau_i\tau}\,.
\label{chj}
\qqq
The characters $\,\chi_i\,$ may be used to define flat line bundles 
$\,\hat L_i\,$ over groups $\,G_i\,$ by setting
\qq
\hat L_i\ =\ (\tilde G_i\times\NC)/\mathop{\sim}\limits_{^i}
\label{Lt}
\qqq
with the equivalence relation
\qq
(\tilde\gamma,u)\ \ \mathop{\sim}\limits_{^i}\ \ (\tilde\gamma\zeta,
\,\chi_i(\zeta)^{-1}u)\,
\label{again}
\qqq
for $\,\zeta\in\CZ_i$. \,The circle subbundle of $\,\hat L_i\,$ forms 
under the multiplication induced by the point-wise one 
in $\,\tilde G_i\times U(1)\,$ a central extension $\,\hat G_i\,$ 
of $\,G_i$. \,These extensions were a centerpiece of the construction
of ref.\,\,\cite{Meinr}. 
\vskip 0.4cm

For $\,I\subset R\,$ with more than one element,
one defines subgroups $\,G_{I}\subset G\,$ as the adjoint action 
stabilizers of elements $\,\ee^{2\pi i\tau}\,$ with $\,\tau\,$ in 
the open simplex in $\,\CA\,$ generated by vertices $\,\tau_i,\ i\in I$
\,($G_I\,$ does not depend on the choice of $\,\tau$). In general, 
$\,G_I\not=\mathop{\cap}\limits_{i\in I} G_i$.
\,To spare on notation, we shall write $\,G_{\{i,j\}}
=G_{ij}\,$ with $\,G_{ii}=G_i\,$ etc.   
\,Let $\Ng_I\,$ be the Lie algebra of $\,G_I\,$ and $\,\tilde G_I\,$ 
its universal cover such that $\,G_I=\tilde G_I/\CZ_I$. \,For $\,J\supset I$, 
$\,G_{J}\subset G_{I}\,$ and the inclusion $\,\Ng_J\subset\Ng_I\,$ 
induces the homomorphisms of the universal covers 
\qq
&&\tilde G_J\ \ \longrightarrow\ \ \tilde G_I \cr
&&\,\downarrow\hspace{1.7cm}\downarrow\label{emb}\\
&& G_J\m\ \ \ \subset\m\ \ \ G_I
\nonumber
\qqq
which map $\,\CZ_J\,$ in $\,\CZ_I$. \,$G_R\,$ is equal to the Cartan 
subgroup $\,T\,$ of $\,G\,$ so that $\,\tilde G_R=\Nt\,$ and for each 
$\,I\,$ one has a natural homomorphism 
\qq
\Nt\ \mathop{\longrightarrow}\limits^{\ee_I^{2\pi i\,\cdot\,}}\ \ \tilde G_I
\label{niso}
\qqq
that maps onto a commutative subgroup $\,\tilde T_I\,$ covering 
$\,T\subset G_I\,$ and sends the coroot lattice $\,Q\,$ onto
$\,\CZ_I$. \,Let
\qq
a_{ij}\ =\ i\,\tr\,(\tau_j-\tau_i)\m(\gamma^{-1}d\gamma)\,.
\label{aij}
\qqq
be a one form on $\,G_{ij}$. \,It is easy to see that $\,a_{ij}\,$
is closed. Indeed,
\qq
da_{ij}\ =\ i\,\tr\,(\tau_i-\tau_j)\m(\gamma^{-1}d\gamma)^2\ =\ 0\,,
\qqq
where the last equality follows from the easy to check fact that 
the adjoint action of the Lie algebra $\,\Ng_{ij}\,$ (and, hence, also
of $\,G_{ij}$) preserves $\,\tau_i-\tau_j$. \,Let $\,\chi_{ij}\,$
be a $\,U(1)$-valued function on the covering group $\,\tilde G_{ij}\,$ such
that $\,i\m\chi_{ij}^{-1}d\chi_{ij}\,$ is the pullback of $\,a_{ij}\,$ 
to $\,\tilde G_{ij}\,$ and that $\,\chi_{ij}(1)=1$. \,Explicitly,
\qq
\chi_{ij}(\tilde\gamma)\ =\ 
\exp\Big[{_1\over^i}\int\limits_{\tilde\gamma}a_{ij}\Big]\,,
\label{rel0}
\qqq
where $\,\tilde\gamma\,$ is interpreted as a homotopy class of paths in
$\,G_{ij}\,$ starting from $\,1$. \,It is easy to see that  $\,\chi_{ij}\,$
defines a 1-dimensional representation of $\,\tilde G_{ij}$:
\qq
\chi_{ij}(\tilde\gamma{\tilde\gamma}')\ =\ 
\chi_{ij}(\tilde\gamma)
\,\chi_{ij}(\tilde\gamma')
\label{rel1}
\qqq
and that for $\,\tilde\gamma\in\tilde G_{ijk}\,$ that may be also viewed
as an element of $\,\tilde G_{ij},\ \,\tilde G_{jk}\,$ and $\,\tilde G_{ik}\m$,
see diagram (\ref{emb}),   
\qq
\chi_{ij}(\tilde\gamma)\,\chi_{jk}(\tilde\gamma)\ =\ 
\chi_{ik}(\tilde\gamma)\,.
\label{rel3}
\qqq
As may be easily seen from the definition (\ref{rel0}), 
\qq
\chi_{ij}(\ee_{ij}^{2\pi i\tau})\ =\ \ee^{\m2\pi i\,
\tr\,(\tau_j-\tau_i)\m\tau}
\label{chij}
\qqq
for $\,\tau\in\Nt$. \,In particular, for $\,\zeta\in \CZ_{ij}\,$ also treated 
as an element of $\,\CZ_i\,$ and $\,\CZ_j$, 
\qq
\chi_{ij}(\zeta)\ =\ \chi_j(\zeta)\,\chi_i(\zeta)^{-1}\,.
\label{rel2}
\qqq
\vskip 0.3cm

The construction of the basic gerbe $\,\CG=(Y,B,L,\mu)\,$ over group $\,G\,$ 
described in \cite{Meinr} uses a specific open covering $\,(\CO_i)\,$ of $\,G$, 
\,where
\qq
\CO_i\ =\ \{\, h\,\ee^{2\pi i\tau}{h}^{-1}\,\,|\,\,h\in G\,,\ \,
\tau\in\CA_i\,\}\,.
\qqq
Over sets $\,\CO_i\m$, \m the closed 3-form $\,H\,$ becomes exact. 
More concretely, the formulae
\qq
B_i\ =\ {_1\over^{4\pi}}\,\tr\,(h^{-1}dh)\,\ee^{2\pi i\tau}\,
(h^{-1}dh)\,\ee^{-2\pi i\tau}\ +\ i\,\tr\,(\tau-\tau_i)(h^{-1}dh)^2
\label{Bi}
\qqq
define smooth 2-forms on $\,\CO_i\,$ such that $\,dB_i=H$.
\,More generally, let $\,\CO_I=\mathop{\cap}\limits_{i\in I}\m\CO_i$.
\,Since the elements $\,\ee^{2\pi i\tau}\,$ with $\,\tau\in\CA_I\,$ 
have the adjoint action stabilizers contained in $\,G_I$, \,the maps
\qq
\CO_I\,\ni\,g\,=\,h\,\ee^{2\pi i\tau}\m h^{-1}\ \mathop{\longrightarrow}
\limits^{\rho_I}\ \ hG_I\,\in\,G/G_I
\label{rhos}
\qqq
are well defined. They are smooth \cite{Meinr}. They will play an important 
role below. \,On the double intersections $\,\CO_{ij}\m$,
\qq
B_j-B_i\ =\ i\,\tr\,(\tau_i-\tau_j)(h^{-1}dh)^2
\qqq
are closed 2-forms but, unlike in the case of the $\,SU(N)\,$ (and $\,Sp(2N)$) groups,
their periods are not in $\,2\pi\NZ,$ in general.  As a result, they are
not curvatures of line bundles over $\,\CO_{ij}$. \,It is here that the
general case departs from the $\,SU(N)\,$ one as described in \cite{GR} where 
it was enough to take $\,Y=\sqcup\CO_i$. \,Instead, following \cite{Meinr},
we consider manifolds 
\qq
Y_i\ =\ \{\,(g,h)\,\in\,\CO_i\times G\ |\ \rho_i(g)\,=\,hG_i\,\}\,.
\qqq
With the natural projections $\,\pi_i\,$ on $\,\CO_i\,$ and
the right actions of $\,G_i\,$ on \,$\,h_i$, $\,Y_i\,$ are principal 
$\,G_i$-bundles over $\,\CO_i$. \,We set
\qq
Y\ \,=\ \,\mathop{\sqcup}\limits_{i=0,..\m,r}Y_i\,
\qqq 
with the projection $\,\pi:Y\rightarrow G\,$ that restricts
to $\,\pi_i\,$ on each $\,Y_i$. \,The curving 2-form $\,B\,$ on $\,Y\,$ is 
defined by setting
\qq
B|_{_{Y_i}}\ =\ \pi_i^*B_i\,.
\qqq
Clearly, $\,dB=\pi^*H\,$ \,as required.
\vskip 0.4cm

We are left with the construction of the line bundle with connection $\,L\,$
over $\,Y^{[2]}\,$ and of the groupoid product $\,\mu$.
\,Note that
\qq
Y^{[2]}\ =\ \mathop{\sqcup}\limits_{(i,j)\atop i,j=0,..\m,r}Y_{ij}\,,
\ \quad{\rm where}\ \ 
Y_{ij}\ =\ Y_i\times_{_{\CO_{ij}}}\hspace{-0.14cm}Y_j\,.
\qqq
It will be convenient to use another description of $\,Y_{ij}$.
\,We may identify
\qq
Y_{ij}\ \ \equiv\ \ \hat Y_{ij}/G_{ij}\,,
\qqq
where 
\qq
\hat Y_{ij}\ =\ \{\,(g,h,\gamma,\gamma')
\,\in\,\CO_{ij}\times
G\times G_i\times G_j\ |\ \rho_{ij}(g)\,=\,hG_{ij}\,\}
\label{ypr}
\qqq
and $\,G_{ij}\,$ acts on $\,\hat Y_{ij}\,$ by the simultaneous right 
multiplication of $\,h,\ \gamma\,$ and $\,\gamma'$. 
\,The identification comes from writing 
\qq
y\ =\ (g,h\gamma^{-1})\,,
\qquad y'\ =\ (g,h{\gamma'}^{-1})
\qqq
for $\,(y,y')\in Y_{ij}$. \,It also works for $\,i=j$.
\,Let us denote by $\,\hat p,\ \hat p_i\,$
and $\,\hat p_j\,$ the projections from $\,\hat Y_{ij}\,$ to
$\,G$, $\,G_i\,$ and $\,G_j$, \,respectively. \,Consider on $\,\hat Y_{ij}\,$
the line bundle
\qq
\hat L_{ij}\ =\ {\hat p}^{\m*}\hat L\m\otimes\,{\hat p_i}^{\m*}{\hat L_i}^{^{-1}}\otimes
\,{\hat p_j}^{\m*}{\hat L}_j\,,
\qqq
where 
$\,\hat L=G\times\NC\,$ is the trivial line bundle over 
$\,G\,$ with the connection form
\qq
A_{ij}\ =\ i\,\tr\,(\tau_j-\tau_i)\m(h^{-1}dh)
\label{Aij}
\qqq
and $\,\hat L_i$, $\,\hat L_j\,$ are the flat line bundles defined in eq\,\,(\ref{Lt}). 
Explicitly, the elements of $\,\hat L_{ij}\,$ may be represented
by the triples $\,(g,h,[\tilde\gamma,{\tilde\gamma}',u]_{_{ij}})\,$
with the equivalence classes corresponding to the relation
\qq
(\tilde\gamma,{\tilde\gamma}',\m u)
\ \ \mathop{\sim}_{^{ij}}\ \ (\tilde\gamma\zeta,
{\tilde\gamma}'\zeta',\m\chi_i(\zeta)\,\chi_j(\zeta')^{-1} u)
\qqq
for $\,\tilde\gamma_i\in\tilde G_i$, $\,\tilde\gamma_j\in\tilde G_j$,
$\,u\in\NC$, $\,\zeta\in \CZ_i\,$ and $\,\zeta'\in \CZ_j$. 
\vskip 0.4cm

We shall lift the action
of $\,G_{ij}\,$ on $\,\hat Y_{ij}\,$ to the action on $\,\hat L_{ij}\,$
by automorphisms of the line bundle preserving the connection and
shall set
\qq
L|_{_{Y_{ij}}}\ =\ \hat L_{ij}/G_{ij}\,\equiv\,L_{ij}\,.
\label{qb}
\qqq
The action of $\,G_{ij}\,$ on $\,\hat L_{ij}\,$ is defined as follows.
\,Let $\,{\tilde\gamma}''\,$ 
be an element in $\,\tilde G_{ij}\,$ projecting to $\,\gamma''\in 
G_{ij}$. \,Consider the map
\qq
(h\m,\,[\tilde\gamma,{\tilde\gamma}',u]_{_{ij}})\ \ \longmapsto\ 
\ (h\gamma''\m,\,[\tilde\gamma{\tilde\gamma}'',{\tilde\gamma}'
{\tilde\gamma}'',\m\chi_{ij}({\tilde\gamma}'')^{-1} u]_{_{ij}})\,,
\label{map}
\qqq
where $\,\chi_{ij}\,$ is given by eq.\,\,\,(\ref{rel0}).
Due to relations (\ref{rel1}) and (\ref{rel2}), the equivalence class 
on the right hand side depends only on $\,\gamma''$. \,Maps (\ref{map})
define the right action of $\,G_{ij}\,$ on $\,\hat L_{ij}\,$ 
(acting trivially on the $\,\CO_{ij}\,$ component). The relation 
between 1-form $\,A_{ij}\,$ and $\,\chi_{ij}\,$ implies that the connection on
$\,\tilde L_{ij}\,$ is preserved by that action. We may then define the
quotient line bundle $\,L_{ij}\,$ by (\ref{qb}). Note that the curvature of
$\,L_{ij}\,$ is given by the closed 2-form $\,F_{ij}\,$ on $\,Y_{ij}\,$
that pulled back to $\,\hat Y_{ij}\,$ becomes
\qq
\hat F_{ij}\ =\ i\,\tr\,(\tau_i-\tau_j)\m(h^{-1}dh)^2\,.
\label{2for}
\qqq
Let $\,p_i\,$ and $\,p_j\,$ denote the natural projections 
of $\,Y_{ij}\,$ on $\,Y_i\,$ and $\,Y_j$, respectively. The required relation 
\qq
{p_j}^*\pi_j^*B_j\ -\ {p_i}^*\pi_i^*B_i\ =\ F_{ij}
\qqq
between the curving 2-form and the curvature of $\,L_{ij}\,$
follows from the comparison of eqs.\,\,(\ref{Bi}) and (\ref{2for})
with the use of the relation $\,\rho_{ij}(g)=hG_{ij}$, \,see (\ref{ypr}).
For $i=j$, line bundle $\,L_{ij}\,$ is flat.
\vskip 0.4cm

We still have to define the groupoid product $\,\mu\,$ in the line bundles 
over 
\qq
Y^{[3]}\ =\ \mathop{\sqcup}\limits_{i,j,k=0,..\m,r}Y_{ijk}\,,
\ \quad{\rm where}\ \ 
Y_{ijk}\ =\ Y_i\times_{_{\CO_{ijk}}}\hspace{-0.14cm}Y_j\times_{_{\CO_{ijk}}}
\hspace{-0.14cm}Y_k\,.
\qqq
For $\,(y,y',y'')\in Y_{ijk}$,
\qq
y\ =\ (g,h\gamma^{-1})\,,\qquad y'\ =\ (g,h{\gamma'}^{-1})\,,
\qquad y''=(g,h{\gamma''}^{-1})\,, 
\label{tbul}
\qqq
where $\,g\in\CO_{ijk}$, $\,h\in G\,$
with $\,\rho_{ijk}(g)=hG_{ijk}\m$, \m and where $\,\gamma\in G_i$, 
$\,\gamma'\in G_j$, $\,\gamma''\in G_k$. 
\,The elements in the corresponding fibers of $\,L_{ij}$, $\,L_{jk}\,$ 
and $\,L_{ik}\,$
may be defined now as the $\,G_{ijk}$-orbits since $\,h\,$
is defined by $\,g\in\CO_{ijk}\,$ up to right multiplication by elements of 
$\,G_{ijk}$. \,Let
\qq
&\displaystyle{\ell_{ij}\ =\ (h\m,\,[{\tilde\gamma}\m,\m{\tilde
\gamma}',\m u]_{_{ij}})G_{ijk}\,\in\,L_{ij}\,,\qquad
\ell_{jk}\ =\ (h\m,\,[{\tilde\gamma}'\m,\m{\tilde\gamma}'',\m u']_{_{jk}})
G_{ijk}\,\in\,L_{jk}\,.}&\cr\cr
&\displaystyle{\ell_{jk}\ =\ (h\m,\m[{\tilde\gamma}\m,\m{\tilde\gamma}'',
\m u\m u']_{_{ijk}})G_{ijk}\,\in\,L_{ik}\,.}&
\nonumber
\qqq
One sets
\qq
\mu\Big(\ell_{ij}\otimes\m\ell_{jk}\Big)\ =\ \ell_{jk}\,.
\qqq
It is easy to see that the right hand side is well defined. 
Checking that $\,\mu\,$ preserves the connection and 
is associative over $\,Y^{[4]}=\mathop{\sqcup}\limits_{i,j,k,l}Y_{ijkl}\,$
is also straightforward (the latter is done
by rewriting the line bundle elements as $\,G_{ijkl}$-orbits). 
\vskip 0.4cm

For $\,\sk\in\NZ$, \,the powers $\,\CG^\sk\,$ of the basic gerbe 
may be constructed the same way by simply exchanging the characters 
$\,\chi_{i}\,$ and homomorphisms $\,\chi_{{ik}}\,$ by their 
$\,\sk^{\m\rm th}\,$ powers and by multiplying the connection forms, 
curvings and curvatures by $\,\sk$. \,Below, we shall use the notation 
$\,[\cdots]^\sk_{_i}\,$ and $\,[\cdots]_{_{ij}}^\sk\,$ for the 
corresponding equivalence classes with such modifications.

\nsection{Basic gerbe on compact non-simply connected groups} 

Let $\,G\,$ be as before and let $\,Z\,$ be a (non-trivial) subgroup of its center.
Let $\,H'\,$ be the 3-form on the non-simply connected group $\,G'=G/Z\,$ that 
pulls back to the 3-form $\,H\,$ on $\,G$. \,We shall construct in this section 
the basic gerbe $\,\CG'=(Y',B',L',\mu')\,$ over group $\,G'\,$ 
with curvature $\,\sk H'$, \,where the level $\,\sk\,$ takes the lowest positive 
(integer) value.
\vskip 0.4cm

Group $\,Z\,$ acts on the Weyl alcove $\,\CA\,$ in the Cartan algebra
of $\,G\,$ by affine transformations. The action is induced from that
on $\,G\,$ that maps conjugacy classes to conjugacy classes and it may be
defined by the formula 
\qq
z\,\ee^{2\pi i\tau}\ =\ w_z^{-1}\ee^{2\pi i\m z\tau}w_z
\qqq
for $\,z\in Z\,$ and $\,w_z\,$ in the normalizer $\,N(T)\subset G\,$
of the Cartan subgroup $\,T$. \,In particular, $\,z\m\tau_i=\tau_{zi}\,$ 
for some permutation $\,i\mapsto zi\,$ of the set $\,R=\{0,1,..\m, r\}\,$
that induces a symmetry of the extended Dynkin diagram with vertices 
belonging to $\,R\,$ and $\,k_{zi}=k_i$, $\,k^{^{\hspace{-0.03cm}\vee}}_{zi}
=k^{^{\hspace{-0.03cm}\vee}}_{i}$. \,Explicitly,
\qq
z\tau\ =\ w_z\m\tau\, w_z^{-1}\,+\,\tau_{z0}\,.
\label{nla}
\qqq
Elements $\,w_z\in N(T)\,$ are defined up to multiplication 
(from the right or from the left) by elements in $\,T,$ so that
their classes $\,\omega_z\,$ in the Weyl group $\,W=N(T)/T\,$ are
uniquely defined. The assignment $\,Z\ni z\,\mathop{\longmapsto}
\limits^{\omega}\,
\omega_z\in W\,$ is an injective homomorphism. However, one cannot always 
choose $\,w_z\in N(T)\,$ so that $\,w_{zz'}=w_{z}w_{z'}$. 
\,The $\,T$-valued discrepancy
\qq
c_{z,z'}\ =\ w_z w_{z'}\m w_{zz'}^{-1}
\label{2coc}
\qqq
satisfies the cocycle condition
\qq
(\delta c)_{z,z',z''}\ \equiv\ (w_z\m c_{z',z''}w_z^{-1})\, c_{zz',z''}^{\,-1}\,
c_{z,z'z''}\,c_{z,z'}^{,-1}\ =\ 1
\label{bcoc}
\qqq
and defines a cohomology class $\,[c]\in H^2(Z,T)\,$ that 
is the obstruction to the existence of a multiplicative choice of $\,w_z\,$
(for a quick r\'esum\'e  of discrete group cohomology, see Appendix A of
\cite{GR}). \,Class $\,[c]\,$ is the restriction to $\,\omega(Z)\subset W\,$ 
of the cohomology class in $\,H^2(W,T)\,$ that characterizes 
up to isomorphisms the extension
\qq
1\ \longrightarrow\ T\ \longrightarrow\ N(T)\ \longrightarrow\ 
W\ \longrightarrow\ 1
\label{cex}
\qqq
that was studied in ref.\,\,\cite{Tits}. The results of \cite{Tits} could
be used to find the 2-cocycle whose cohomology class characterizes the
extension (ref{cex}) and then, by restriction, to calculate $\,c$.
In practice, we found it simpler to obtain the 2-cocycle $\,c\,$ directly,
see Sect.\,\,4.  
\vskip 0.3cm

Let us choose elements $\,e_{z,z'}\in \Nt\,$ such that $\,c_{z,z'}
=\ee^{2\pi i\,e_{z,z'}}$ (of course, they are fixed modulo the
coroot lattice $\,Q^{^{\hspace{0.03cm}\vee}}$). Then
\qq
(\delta e)_{z,z',z''}\ \equiv\ (w_z\m e_{z',z''}w_z^{-1})\m
-\m e_{zz',z''}+\m e_{z,z'z''}
-e_{z,z'}
\label{Bock}
\qqq
is a 3-cocycle on $\,Z\,$ with values in $\,Q^{^{\hspace{0.03cm}\vee}}$. 
\,It defines a cohomology class 
$\,[\delta e]\in H^3(Z,Q^{^{\hspace{0.03cm}\vee}})$, \,the Bockstein image
of $\,[c]\,$ induced by the exact sequence
\qq
0\ \longrightarrow\ Q^{^{\hspace{0.03cm}\vee}}\ \longrightarrow\ \,\Nt\ \,
\mathop{\longrightarrow}\limits^{\ee^{2\pi i\,\cdot\,}}\ \,T\ \m
\longrightarrow\ 1\,.
\qqq
Below, we shall employ for $\,I\subset R\,$ the lifts
\qq
\tilde c_{z,z'}\ =\ \ee_I^{2\pi i\m e_{z,z'}}\ \in\ \tilde T_I\,\subset\,
\tilde G_I
\label{tczz}
\qqq
of $\,c_{z,z'}\,$ to the subgroups $\,\tilde T_I$, \,see (\ref{niso}).
Note that
\qq
(\delta\tilde c)_{z,z',z''}\,\equiv\,\ee_I^{2\pi i\m(\delta e)_{z,z',z''}}
\label{wib}
\qqq
belongs to $\,\CZ_I\subset \tilde T_I$.
\vskip 0.4cm

The structures introduced in the preceding section behave naturally 
under the action of $\,Z$. \,We have
\qq
z\CA_I\ =\ \CA_{zI}\,,\qquad z\CO_{I}\ =\ \CO_{zI}\,,\qquad w_z\m 
G_{I}\m w_z^{-1}\ =\ G_{zI}\,.
\qqq 
The adjoint action of $\,w_z\,$ maps also $\,\Ng_I\,$ onto $\,\Ng_{zI}\,$
and hence lifts to an isomorphism from $\,\tilde G_I\,$ to $\,\tilde G_{zI}\,$
that maps $\,\CZ_I\,$ onto $\,\CZ_{zI}\,$ and for which we shall still use the
notation $\,\tilde\gamma_I\mapsto w_z\tilde\gamma_Iw_z^{-1}$.
\,The maps $\,\CO_i\ni g\mapsto zg\in\CO_{zi}\,$ satisfy the relation
\qq
z^*B_{zi}\ =\ B_i\,.
\label{trb}
\qqq
They may be lifted to the maps
\qq
Y_i\ni y=(g,h)\ \ \longmapsto\ \ z\m y=(zg,\m h\m w_z^{-1})\in Y_{zi}
\qqq 
of the principal bundles $\,Y_i$. \,Note that if $\,c_{z,z'}\not=1\,$ then the lifts 
do not compose.
\vskip 0.4cm

Proceeding to construct the basic gerbe $\,\CG'=(Y',B',L',\mu')\,$ over group $\,G'$, 
\,we shall set 
\qq
Y'\ =\ Y\ =\ \mathop{\sqcup}\limits_{i=0,..\m,r}Y_i\,,\qquad 
B'|_{_{Y_i}}=\sk\,\pi_i^*B_i\,,
\qqq
where $\,Y'\,$ is taken with the natural projection 
$\,\pi'\,$ on $\,G'$. \,Note that a sequence 
$\,(y,y',\dots,y^{(n-1)})\,$  belongs to $\,{Y'}^{[n]}\,$ if
$\,\pi(y)=z\pi(y')=\,...=zz'\cdots z^{(n-2)}\pi(y^{(n-1)})\,$ for
some $\,z,z',\dots,z^{(n-2)}\in Z$. \,Then
\qq
(y,zy',\dots,z(z'(\cdots(z^{(n-2)}y^{(n-1)})\cdots)))\ \ \in\ \ Y^{[n]} 
\qqq
and we may identify
\qq
{Y'}^{[n]}\ \ \cong\ \ \mathop{\sqcup}\limits_{(z,z',\dots,z^{(n-1)})\,\in\,
Z^{n-1}}\,Y^{[n]}\,.
\label{id}
\qqq
Let $\,L'\,$ be the line bundle on $\,{Y'}^{[2]}\,$
that restricts to $\,L^{^\sk}\,$ on each component $\,Y^{[2]}\,$ in the 
identification (\ref{id}), i.e.\,\,to $\,L_{ij}^{^\sk}\,$ on $\,Y_{ij}\subset Y^{[2]}$.
It is easy to see with the use of eq.\,\,(\ref{trb}) that the curvature 
$\,F'\,$ of $\,L'\,$ satisfies the required relation
\qq
F'\ =
\ {p'_2}^*B'\ -\ {p'_1}^*B'\,
\qqq
where, as usual, $\,p'_1\,$ and $\,p'_2\,$ are projections 
in $\,{Y'}^{[2]}\,$ on the first and the second factor. 
\vskip 0.4cm

It remains to define the groupoid multiplication $\,\mu'$.
\,Let $\,(y,y',y'')\in{Y'}^{[3]}\,$ be such that 
$\,(y,zy',z(z'y''))\in Y_{ijk}\subset Y^{[3]}$.
\,We may then write, see (\ref{tbul}),
\qq
y\ =\ (g,\m h\m\gamma^{-1})\,,\quad\ z\m y'\ =\ 
(g,\m h\m{\gamma'}^{-1})\,,\quad\ 
z(z'y'')\ =\ (g,\m h\m{\gamma''}^{-1})
\label{wmr0}
\qqq
with $\,g\in\CO_{ijk}$, $\,h\in G\,$ such that 
$\,\rho_{ijk}(g)=hG_{ijk}\,$ and with $\,\gamma\in G_i$, $\,\gamma'\in G_j$, 
$\,\gamma''\in G_k$. 
\,We shall use the notation $\,i_z\equiv z^{-1}i$, 
$\,\gamma_z\equiv w_z^{-1}\gamma\,w_z\in G_{i_z}\,$ for $\,\gamma\in G_i\,$ 
and $\,\tilde\gamma_z\equiv w_z^{-1}\tilde\gamma\,w_z\in\tilde G_{i_z}\,$ for
$\,\tilde\gamma\in\tilde G_i$. \,Note that
\qq
&&\hbox to 5.5cm{$\displaystyle{y'\ =\ (z^{-1}g,\m h\,w_z
\m{\gamma'_z}^{-1})\,,}
$\hfill}y''\ =\ ((zz')^{-1}g,h\,w_zw_{z'}(\gamma''_z)_{z'}^{-1})\,,\label{wmr1}\\ \cr
&&\hbox to 5.5cm{$\displaystyle{z'y''\ = 
\ (z^{-1}g,\m h\,w_z\m{\gamma''_z}^{-1})\,,}$\hfill}(zz')y'' 
\ =\ (g,\m h\,(c_{z,z'}^{\,-1}\gamma'')^{-1})\,.
\label{wmr2}
\qqq
Employing the explicit description of the line bundles $\,L_{ij}\,$ 
with $\,\tilde\gamma\in\tilde G_i\,$ projecting to $\,\gamma\,$ etc.,
we take the elements
\qq
&&\hbox to 7.5cm{$\displaystyle{\ell_{ij}\ \ \ \,=\ (g,h,\,[\tilde\gamma,
{\tilde\gamma}',\m u]^{^{\m\sk}}_{_{ij}})G_{ijk}}$\hfill}\in\,\ 
L^{^{\,\sk}}_{(y,zy')}\,\quad\ =\,L'_{(y,y')}\,,\label{ell1}\\ \cr
&&\hbox to 7.5cm{$\displaystyle{\ell_{j_zk_z}\ =
\ (z^{-1}g,\m h\, w_z,\,[{\tilde\gamma}'_{z},
{\tilde\gamma}''_{z}
,\m u']^{^{\m\sk}}_{j_zk_z})G_{i_zj_zk_z}}$\hfill}\in\,\ L^{^{\,\sk}}_{(y',z'y'')}\,\ \ =\,
L'_{(y',y'')}\,,\qquad\label{ell2}\\ \cr
&&\hbox to 7.5cm{$\displaystyle{\ell_{ik}\ \ \ \,=\ (g,h,\,[\tilde\gamma,
{\tilde c}_{z,z'}^{\,-1}{\tilde\gamma}'',\m u\m u']^{^{\m\sk}}_{ik})
G_{ijk}}$\hfill}\in\,\ L^{^{\,\sk}}_{(y,(zz')y'')}\,
=\,L'_{(y,y'')}
\label{ell3}
\qqq
where $\,\tilde c_{z,z'}\in\tilde G_k\,$ is given by eq.\,\,(\ref{tczz}).
Then necessarily, 
\qq
\mu'\m(\m\ell_{ij}\otimes\ell_{j_zk_z})\ =\ u^{ijk}_{z,z'}\,\ell_{ik}\,,
\label{mup}
\qqq
where
%$\,{\tilde c}_{z,z'}\,$ is a lift of $\,c_{z,z'}\in T\,$ to
%the Cartan subgroup of $\,\tilde G_k\,$ and 
$\,u_{ijk}^{z,z'}\,$ are
numbers in $\,U(1)$. 
%\,Note that the simultaneous change
%$\,{\tilde c}_{z,z'}\mapsto{\tilde c}_{z,z'}\zeta\,$ of the lift 
%of $\,c_{z,z'}\,$ 
%with $\,\zeta\in\CZ_k\,$ and $\,u^{ijk}_{z,z'}\mapsto\chi_{k}(\zeta)^{^{\sk}}
%\,u^{ijk}_{z,z'}\,$ does not change the right hand side of (\ref{mup})
%so that it depends only on
%\qq
%c_{z,z'}^{ijk}\ =\ [\tilde c_{z,z'},u^{ijk}_{z.z'}]^{^{-\sk}}_{_{k}}\,\in
%\hat G^{^{-\sk}}_k\subset\hat L^{^{-\sk}}_k\,,
%\qqq
%where $\,\hat G^{^{-\sk}}_k$, the circle bundle of $\,L^{^{-\sk}}_k$, \,is
%the $\,-\sk^{\rm th}\,$ power of the central extension $\,\hat G_k\,$
%of $\,G_k$. 
That the right hand side of the definition (\ref{mup}) does not depend 
on the choice of the representatives of the classes on the left hand side 
follows from 
\vskip 0.7cm

\noindent{\bf Lemma 1.} \ For $\,z\in Z$, $\,\zeta\in\CZ_i\,$ and  
$\,\tilde\gamma\in\tilde G_{ij}$,
\qq
\chi_{i_z}(\zeta_z)\ =\ \chi_i(\zeta)\,,\qquad
\chi_{i_zj_z}(\tilde\gamma_z)\ =\ \chi_{ij}(\tilde\gamma)\,.
\label{toc}
\qqq
\vskip 0.4cm

\noindent{\bf Proof.} \ Let $\,\zeta=\ee_i^{2\pi i p}\,$ for 
$\,p\in P_i^{^{\hspace{-0.03cm}\vee}}\cap Q^{^{\hspace{-0.03cm}\vee}}$.
Then $\,\zeta_z=w_z^{-1}\zeta\m w_z=\ee_i^{2\pi i\m w_z^{-1}p\,w_z}\,$ and
\qq
\chi_{i_z}(\zeta_z)\ =\ \ee^{2\pi i\,\tr\,
w_z^{-1}(\tau_i-\tau_{z0})\, w_z\,w_z^{-1}p\, w_z}
\ =\ \ee^{2\pi i\,\tr\,\tau_{i}p}\ =\ \chi_i(\zeta)\,
\qqq
where we used the fact that $\,\tau_{z0}=\lambda_{z0}$. \,The second
relation in (\ref{toc}) follows immediately from the definition (\ref{rel0})
of $\,\chi_{ij}\,$ and the identity $\,\tau_{j_z}-\tau_{i_z}
=w_z^{-1}(\tau_j-\tau_i)\,w_z$. 
\vskip 0.2cm

\hspace{10cm}$\Box$
\vskip 0.4cm

It remains to find the conditions under which $\,\mu'\,$ is associative.
In Appendix A, we show by an explicit check that associativity of
$\,\mu'\,$ requires that
\qq
u^{j_zk_zl_z}_{z',z''}\,(u^{ikl}_{zz',z''})^{\,-1}\,u^{ijl}_{z,z'z''}\,
\m(u^{ijk}_{z,z'})^{-1}\ =\ \chi_{_{kl}}(\tilde c_{z,z'})^{^{\sk}}\,
\chi_{_l}((\delta\tilde c)_{z,z',z''})^{^{\sk}}\,.
\label{rtc}
\qqq
This provides an extension of the relation (4.6) of \cite{GR}
obtained for $\,G=SU(N)\,$ and it may be treated similarly.
First, we set
\qq
u_{z,z'}\ =\ u_{z,z'}^{(0)\m(z0)\m(zz'0)}
\qqq
and observe that, for $\,i=j_z=k_{zz'}=l_{zz'z''}=0$, \,eq.\,\,(\ref{rtc})
reduces to the cohomological equation   
\qq
\delta u\ =\ U\,,
\label{rtc1}
\qqq
where
\qq
(\delta u)_{z,z',z''}&\equiv& u_{z',z''}\,u_{zz',z''}^{\,-1}
\,u_{z,z'z''}\,u_{z,z'}^{\,-1}
\qqq
is the coboundary of the $\,U(1)$-valued 2-chain on $\,Z\,$ and 
\qq
U_{z,z',z''}\ =\ \chi_{_{(zz'0)\m(zz'z''0)}}(\tilde c_{z,z'})^{^{\sk}}\,
\chi_{_{zz'z''0}}((\delta\tilde c)_{z,z',z''})^{^{\sk}}\,.
\qqq
More exactly, with the use of formulae 
(\ref{chij}), (\ref{chj}), (\ref{tczz}), (\ref{wib}), (\ref{bcoc})
and (\ref{nla}), one obtains:
\qq 
U_{z,z',z''}&=&\ee^{\,2\pi i\,\sk\,\tr\,[\,(\tau_{zz'z''0}-z\tau_{zz'0})
e_{z,z'}\,+\,\tau_{zz'z''0}(w_z e_{z',z''}w_z^{-1}-e_{zz',z''}
+e_{z,z'z''}-e_{z,z'})\,]}\cr
&=&\ee^{\,2\pi i\,\sk\,\tr\,[\,(\tau_{z'z''0}
-\tau_{z^{-1}0})e_{z',z''}\,-\,\tau_{zz'0}e_{z,z'}\,-\,\tau_{zz'z''0}
(e_{zz',z''}-e_{z,z'z''})\,]}\,.
\label{Uzzz}
\qqq 
\noindent The cohomological equation (\ref{rtc1}) is consistent since  
\vskip 0.7cm

\noindent{\bf Lemma 2.} \ $U_{z,z',z''}\,$ defines a $\,U(1)$-valued
3-cocycle on $\,Z\,$:
\qq
(\delta U)_{z,z',z'',z'''}\ \equiv\ 
U_{z',z'',z'''}\,U_{zz',z'',z'''}^{\,-1}\,U_{z,z'z'',z'''}\,
U_{z,z',z''z'''}^{\,-1}\,U_{z,z',z''}\ =\ 1\,.
\label{tcc}
\qqq
\vskip 0.4cm

\noindent It is enough to analyze the condition (\ref{rtc1}) due to 
the following 
\vskip 0.7cm

\noindent{\bf Lemma 3.}\ \ Let $\,(u_{z,z'})\,$ be a solution of 
eq.\,\,(\ref{rtc1}). Then
\qq
u_{z,z'}^{ijk}\ =\ \chi_{_{k\m(zz'0)}}(\tilde c_{z,z'})^{^{-\sk}}\,
u_{z,z'}\,,
\label{tl2}
\qqq
solves eq.\,\,(\ref{rtc}). 
\vskip 0.7cm

\noindent{\bf Remark.\ \ }That the push-forward
of a gerbe by a covering map $\,M\mapsto M/\Gamma\,$ requires solving 
a cohomological problem $\,U=\delta u\,$ for a $\,U(1)$-valued 3-cocycle 
$\,U\,$ on discrete group $\,\Gamma$, with $\,[U]\in H^3(\Gamma,U(1))\,$
describing the obstruction class, is a general fact, see \cite{RT}. 
Similar cohomological equation, but in one degree less, with obstruction
class in $\,H^2(\Gamma,U(1))$, \,describes pushing forward a line bundle. 
As for the relation (\ref{tl2}), it is of a geometric origin, as has been 
explained in \cite{GR}: if we choose naturally a stable isomorphism 
between $\,\CG^{^\sk}\,$ and $\,(z^{-1})^{^*}\CG^{^\sk}\,$ then 
the elements $\,\ell_{_{ij}}^{-1}\otimes\ell_{_{j_zk_z}}^{-1}\otimes
\ell_{_{ik}}\,$ determine flat sections $\,s_{_{ijk}}\,$ of a flat line 
bundle $\,R^{z,z'}\,$ on $\,G$. Sections $\,s_{_{ijk}}\,$ are defined over
sets $\,\CO_{ijk}\,$ and over their intersections $\,\CO_{_{ijk\m i'j'k'}}$,
\,they are related by $\,\,s_{_{i'j'k'}}=\chi_{_{k'k}}(\tilde c_{z,z'})^{^{\sk}}
\m s_{_{ijk}}\m.$
\vskip 0.6cm

Proofs of Lemmas 2 and 3 may be found in Appendix B.
\,The obstruction cohomology class $\,[U]\in H^3(Z,U(1))\,$ 
is $\,\sk$-dependent. The level $\,\sk\,$ of the basic gerbe 
$\,\CG'\,$ over $\,G'\,$ corresponds to the smallest positive value 
for which this 
class is trivial so that eq.\,\,(\ref{rtc1}) has a solution. In the 
latter case, different solutions $\,u\,$ differ by the multiplication 
by a $\,U(1)$-valued 2-cocycle $\,\tilde u$, 
$\,\delta\tilde u=1$. \,If $\,\tilde u\,$ is cohomologically trivial, 
i.e.\,\,$\,\tilde u_{z,z'}=v_{z'}\m v_{zz'}^{-1}\m v_z$, \,then the modified 
solution leads to a stably isomorphic gerbe over $\,G'$. 
\,Whether multiplication of $\,u\,$ by cohomologically non-trivial 
cocycles $\,\tilde u\,$ leads to stably non-isomorphic gerbes depends 
on the cohomology group $\,H^2(G',U(1))\,$ that classifies different 
stable isomorphism classes of gerbes over $\,G'\,$ with fixed curvature. 
This is trivial for all simple groups except for $\,G'=SO(4N)/\NZ_2\,$ 
when $\,H^2(G',U(1))=\NZ_2$, \,see \cite{FGK}.  
\vskip 0.3cm

\nsection{Cocycles $\,c\,$ and $\,U$}
\vskip 0.1cm

It remains to calculate the cocycles $\,c=(c_{z,z'})$, \,elements
$\,e_{z.z'}\in Q^{^{\hspace{-0.03cm}\vee}}\,$ such that $\,c_{z.z'}
=\ee^{2\pi i\m e_{z,z'}}\,$ and the cocycles $\,U=(U_{z,z',z''})$, 
\,see eqs.\,\,(\ref{2coc}) and (\ref{Uzzz}), and to solve 
the cohomological equation (\ref{rtc1}) for all simple, connected,
simply connected groups $\,G\,$ and all subgroups $\,Z\,$ of their 
center. 
\vskip 0.07cm

\subsection{Groups $\,A_r=SU(r+1)\,$}
\vskip 0.05cm

The Lie algebra $\,su(r+1)\,$ is composed of traceless
hermitian $\,(r+1)\times(r+1)$ matrices. The Cartan algebra 
may be taken as the subalgebra of diagonal matrices.
Let $\,e_i$, $i=1,2,\dots,r+1$, \,denote the diagonal matrices
with the $\,j$'s entry $\,\delta_{ij}\,$ with 
$\,\tr\,e_ie_j=\delta_{ij}$. \,Roots and coroots 
of $\,su(r+1)\,$ have then the form $\,e_i-e_j\,$ for $\,i\not=j\,$
%\,The coroot lattice is composed of integer 
%combinations $\,\sum\limits_{i=1}^{r+1}
%n_ie_i\,$ such that $\,\sum\limits_{i=1}^{r+1}n_i=0$.
and the standard choice of simple roots is $\,\alpha_i=e_i-e_{i+1}$.
\,The center is $\,\NZ_{r+1}\,$ and it may be generated by $\,z=\ee^{-2\pi i
\theta}\,$ with $\,\theta=\lambda_r^{\hspace{-0.03cm}^\vee}=
-e_{r+1}+{1\over r+1}\sum\limits_{i=1}^{r+1}e_i$. 
\,The permutation $\,zi=i+1\,$ for $\,i=0,1,\dots,r-1$, $\,zr=0\,$
generates a symmetry of the extended Dynkin diagram:

\leavevmode\epsffile[-85 -20 235 105]{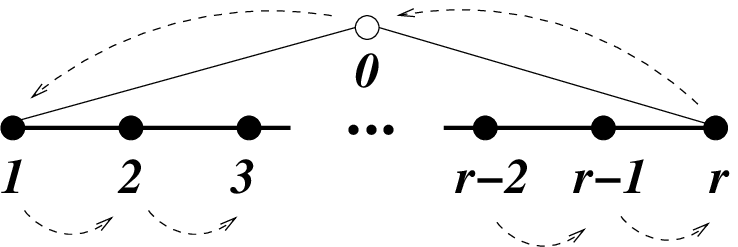}

\noindent The adjoint action of $\,w_z\in N(T)\subset SU(r+1)\,$ on the 
Cartan algebra may be extended to all diagonal matrices
by setting
\qq
w_z\,e_i\,w_z^{-1}\ =\ \cases{\hbox to 1.2cm{$\,e_1$\hfill}{\rm if}\ \ 
i=r+1\,,\cr 
\hbox to 1.2cm{$\,e_{i+1}$\hfill}{\rm otherwise}\,.}
\qqq
It is generated by the product
\qq
\ \ r_{\alpha_1}r_{\alpha_2}\cdots r_{\alpha_r}
\qqq
of $\,r\,$ reflections in simple roots. 
We may take
\qq
\hspace{-0.4cm}w_z\ =\ \ee^{{\pi ir\over r+1}}
\left(\matrix{_0&_0&_0&..&_0&_0&_1\cr_1&_0&_0&..&_0&_0&_0\cr
_0&_1&_0&..&_0&_0&_0\cr.&.&.&..&.&.&.\cr.&.&.&..&.&.&.\cr
_0&_0&_0&..&_1&_0&_0\cr_0&_0&_0&..&_0&_1&_0}\right).
\qqq
Setting $\,w_{z^n}=w_z^n\,$ for $\,n=0,1,\dots,r$, \,we then obtain
\qq
c_{z^n,z^m}\ =\ \cases{\hbox to 1.2cm{$\,1$\hfill}{\rm if}\ \ n+m\leq r\,,\cr
\hbox to 1.2cm{$\,w_z^{r+1}$\hfill}{\rm if}\ \ n+m>r\,.}
\qqq
Since $\,w_z^{r+1}=(-1)^r=\ee^{2\pi iX}\,$ for $\,X={r(r+1)\over2}\theta$, 
\,we may take
\qq
e_{z^n,z^m}\ =\ \cases{\hbox to 1.5cm{$\,0$\hfill}{\rm if}\ \ n+m\leq r\,,\cr
\hbox to 1.5cm{$\,{{r(r+1)}\over2}\theta$\hfill}{\rm if}\ \ n+m>r\,.}
\qqq
Explicit calculation of the right hand side of eq.\,\,(\ref{Uzzz})
gives
\qq
U_{z^n,z^{n'},z^{n''}}\ =\ (-1)^{\,\sk\, r\,n''
{n+n'-[n+n']\over r+1}}\,,
\qqq
where $\,0\leq n,n',n''\leq r\,$ and for an integer $\,m$, 
$\,[m]=m\,\,{\rm mod}\,\,(r+1)\,$ with $\,0\leq [m]
\leq r$. 
\vskip 0.3cm

Let $\,r+1=N'N''$ and $\,Z\,$ be the cyclic subgroup of order $\,N'\,$ 
of the center generated by $\,z^{N''}$. 
\,If $\,N''\,$ is even 
or $\,N'\,$ is odd or $\,\sk\,$ is even, then the restriction to $\,Z\,$
of the cocycle $\,U\,$ is trivial. In the remaining case of $\,N'\,$ even, 
$\,N''\,$ odd and $\,\sk\,$ odd it defines a nontrivial class in 
$\,H^3(Z,U(1))$. \,Hence the smallest positive value of the level
for which the cohomological equation (\ref{rtc1}) nay be solved is
\qq
\sk\ =\ \cases{\,1\qquad{\rm for\ \ }N'\ \ {\rm odd\ \,or}
\ \ N''\ \ {\rm even}\,,\cr
\,2\qquad{\rm for\ \ }N'\ \ {\rm even\ \,and}\ \ N''\ \ {\rm odd}\,,} 
\qqq
in agreement with \cite{GR}. For those values of $\,\sk$, \,one may take 
$\,u_{z^{n},z^{n'}}\equiv1\,$ as the solution of eq.\,\,(\ref{rtc1}). 
\vskip 0.06cm
 
\subsection{Groups $\,B_r=Spin(2r+1)\,$}
\vskip 0.04cm

The Lie algebra of $\,B_r\,$ is $\,so(2r+1)$. \,It is composed of imaginary
antisymmetric $\,(2r+1)\times(2r+1)$ matrices. The Cartan algebra
may be taken as composed of $\,r\,$ blocks 
$\,\left(\matrix{_0&_{-it_i}\cr^{it_i}&^0}\right)\,$ placed diagonally,
with the last diagonal entry vanishing. Let $\,e_i\,$ denote the matrix 
corresponding to $\,t_j=\delta_{ij}$. \,With the invariant form normalized 
so that $\,\tr\,e_ie_j=\delta_{ij}$, \,roots of $\,so(2r+1)\,$ have 
the form $\pm e_i\pm e_j\,$ for $\,i\not=j\,$ and $\pm e_i\,$
and one may choose $\,\alpha_i=e_i-e_{i+1}\,$ 
for $\,i=1,\dots,r-1\,$ and $\,\alpha_r=e_r\,$ as the simple roots.
\,The coroots are $\,\pm e_i\pm e_j\,$ for $\,i\not=j\,$ and $\,\pm 2e_i$. 
%\,The coroot lattice is composed of the elements 
%$\,\sum\limits_{i=1}^rn_ie_i\,$ for integer $\,n_i\,$ such that their sum 
%is even. 
The center of $\,Spin(2r+1)\,$ is $\,\NZ_2\,$
with the non-unit element $\,z=\ee^{-2\pi i\theta}\,$ with $\,\theta=
\lambda_1^{\hspace{-0.03cm}^\vee}=e_1$. $\,SO(2r+1)=Spin(2r+1)/\{1,z\}$. 
The permutation $\,z0=1,\,z1=0,\,zi=i\,$ for $\,i=2,\dots,r\,$
generates a symmetry of the extended Dynkin diagram:

\leavevmode\epsffile[-75 -25 245 90]{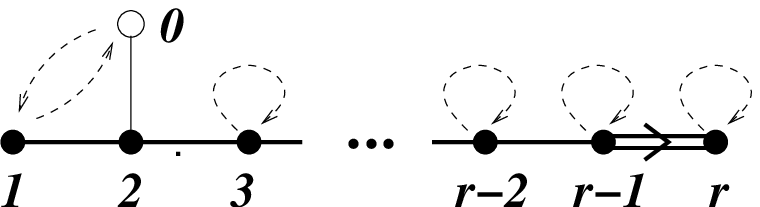}

\noindent The adjoint action of $\,w_z\in N(T)\,$ is given by
\qq
w_z\,e_i\,w_z^{-1}\ =\ \cases{\hbox to 1cm{$\,-e_1$\hfill}{\rm if}\ 
\ i=1\,,\cr
\hbox to 1cm{$\,\,e_i$\hfill}{\rm if}\ \ i\not=1\,.}
\qqq
It may be generated by the product
\qq
r_{\alpha_1}r_{\alpha_2}\cdots r_{\alpha_{r-2}}r_{\alpha_{r-1}}r_{\alpha_r}
r_{\alpha_{r-1}}\cdots r_{\alpha_2}r_{\alpha_1}
\qqq
of $\,2r-1\,$ reflections in simple roots.  
\,Element $\,w_z\,$ may be taken as the lift to $\,Spin(2r+1)$ of the matrix
\qq
\left(\matrix{_1&_0&_0&..&_0&_0\cr_0&_{-1}&_0&..&_0&_0\cr_0&_0&_{-1}&..&_0&_0\cr
.&.&.&..&.&.\cr.&.&.&..&.&.\cr_0&_0&_0&..&_{-1}&_0\cr_0&_0&_0&..&_0&_{-1}}
\right)
\qqq
in $\,SO(2r+1)$. \,Setting also $\,w_1=1$, \,we infer that 
\qq
c_{1,1}\ =\ c_{1,z}\ =\ c_{z,1}\ =\ 1\,,\qquad c_{z,z}\ =\ w_z^2\,.
\qqq
Since $\,w_z^2\,$ projects to $\,1\,$ in $\,SO(2r+1)$, \,it is equal to 
$\,1\,$ or to $\,z$. \,To decide which is the case, we write $\,w_z=
\CO\,\ee^{2\pi iX}\CO^{-1}$, \,where
$\,\CO\in Spin(2r)\,$ projects to the matrix
\qq
\left(\matrix{_0&_0&_0&..&_0&_0&_1\cr_0&_1&_0&..&_0&_0&_0\cr
_0&_0&_1&..&_0&_0&_0\cr.&.&.&..&.&.&.\cr.&.&.&..&.&.&.\cr
%_0&_0&_0&..&_1&_0&_0\cr
_0&_0&_0&..&_0&_1&_0\cr
_{-1}&_0&_0&..&_0&_0&_0}\right)
\qqq
in $\,SO(2r+1)\,$ and $\,\,X={1\over2}\sum\limits_{i=1}^r e_r\,\,$
so that $\,\ee^{2\pi iX}\,$ projects to the matrix
\qq
\left(\matrix{_{-1}&_0&_0&..&_0&_0\cr_0&_{-1}&_0&..&_0&_0\cr
_0&_0&_{-1}&..&_0&_0\cr.&.&.&..&.&.\cr.&.&.&..&.&.\cr
_0&_0&_0&..&_{-1}&_0\cr_0&_0&_0&..&_0&_1}\right)
\qqq
in $\,SO(2r+1)$. \,Now $\,w_z^2=1\,$ if and only if $\,2X\,$ is in
the coroot lattice. This happens if $\,r\,$ is even. We may then take
\qq
e_{1,1}\ =\ e_{1,z}\ =\ e_{z,1}\ =\ e_{z,z}\ =\ 0
\qqq
for even $\,r\,$ and             
\qq
e_{1,1}\ =\ e_{1,z}\ =\ e_{z,1}\ =\ 0\,,\qquad e_{z,z}\ =\ \theta
\qqq
for odd $\,r$. \,\,Here $\,\,U_{z^n,z^{n'},z^{n''}}\,\equiv1\,\,$ for 
all $\,\,0\leq n,n',n''\leq 1$. \,\,Hence $\,\sk=1\,$ and 
$\,\,u_{z^{n},z^{n'}}\equiv 1\,\,$ solves eq.\,\,(\ref{rtc1}).
\vskip 0.15cm

\subsection{Groups $\,C_r=Sp(2r)\,$}
\vskip 0.13cm

This is a group composed of unitary $\m(2r)\times(2r)$ matrices $\,U\,$
such that $\,U^T\Omega U=\Omega\,$ for $\,\Omega\,$ built
of $\,r\,$ blocks $\,\left(\matrix{_0&_1\cr ^{-1}&^0}\right)\,$ placed 
diagonally.
%\qq
%\Omega\ =\ \left(\matrix{_0&_1&_0&_0&..&_0&_0\cr_{-1}&_0&_0&_0&..&_0&_0\cr
%_0&_0&_0&_1&..&_0&_0\cr_0&_0&_{-1}&_0&..&_0&_0\cr.&.&.&.&..&.&.\cr
%.&.&.&.&..&.&.\cr
%_0&_0&_0&_0&..&_0&_1\cr_0&_0&_0&_0&..&_{-1}&_0}\right).
%\qqq
Its Lie algebra $\,sp(2r)\,$ is composed of hermitian $\,(2r)\times(2r)$ 
matrices $\,X\,$ such that $\,\Omega X\,$ is symmetric. The Cartan
subalgebra may be taken as composed of $\,r\,$ blocks 
$\,\left(\matrix{_0&_{-it_i}\cr^{it_i}&^0}\right)\,$ placed diagonally.
Let $\,e_i\,$ denote the matrix 
corresponding to $\,t_j=\delta_{ij}$. \,With the invariant form normalized 
so that $\,\tr\,e_ie_j=2\delta_{ij}$, \,roots of $\,sp(2r)\,$ have 
the form ${1\over2}(\pm e_i\pm e_j)\,$ for $\,i\not=j\,$ and $\pm e_i$.
The simple roots may be chosen as  $\,\alpha_i={1\over 2}(e_i-e_{i+1})\,$ 
for $\,i=1,\dots,r-1\,$ and $\,\alpha_r=e_r$.
\,The coroots are $\,\pm e_i\pm e_j\,$ for $\,i\not=j\,$ and $\,\pm e_i$. 
%\,The coroot lattice is composed of the elements 
%$\,\sum\limits_{i=1}^rn_ie_i\,$ for integer $\,n_i$. 
\,The center of $\,Sp(2r)\,$ is $\,\NZ_2\,$ with the non-unit
element $\,z=\ee^{-2\pi i\theta}\,$ for 
$\,\theta=\lambda_r^{\hspace{-0.03cm}^\vee}={1\over2}
\sum\limits_{i=1}^re_i$. \,The permutation $\,zi=r-i\,$ for $\,i=0,1,
\dots,r\,$ generates a symmetry of the extended Dynkin diagram:

\leavevmode\epsffile[-60 -35 270 130]{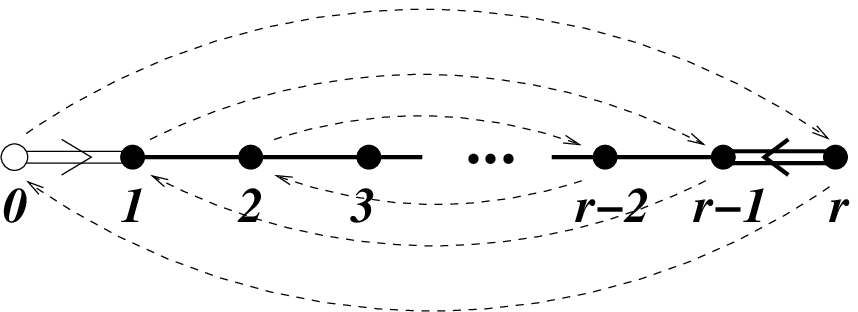}

\noindent Group $\,Sp(2r)\,$ is simply connected. The adjoint action
of $\,w_z\,$ on the Cartan algebra is given by
\qq
w_z\,e_i\,w_z^{-1}\ =\ -e_{r-i+1}\,.
\qqq
It may be generated by the product
\qq
r_{\alpha_r}r_{\alpha_{r-1}}r_{\alpha_r}\cdots r_{\alpha_2}\cdots
r_{\alpha_{r-1}}r_{\alpha_r}
r_{\alpha_1}\cdots r_{\alpha_{r-1}}r_{\alpha_r}
\qqq
of $\,{r(r+1)\over2}\,$ reflections in simple roots. 
Element $\,w_z\,$ may be taken as the matrix
\qq
\hspace{0.3cm}
\left(\matrix{_0&_0&_0&..&_0&_0&_i\cr_0&_0&_0&..&_0&_i&_0\cr
_0&_0&_0&..&_i&_0&_0\cr.&.&.&..&.&.&.\cr.&.&.&..&.&.&.\cr
%_0&_0&_i&..&_0&_0&_0\cr
_0&_i&_0&..&_0&_0&_0\cr_i&_0&_0&..&_0&_0&_0}
\right)
\qqq
in $\,Sp(2r)$. \,Setting also $\,w_1=1$, \,we infer that 
\qq
c_{1,1}\ =\ c_{1,z}\ =\ c_{z,1}\ =\ 1\,,\qquad c_{z,z}\ =\ w_z^2\,=\,-1\,
=\,z
\qqq
so that we may take 
\qq
e_{1,1}\ =\ e_{1,z}\ =\ e_{z,1}\ =\ 0\,,\qquad e_{z,z}\ =\ \theta
\qqq
which results in
\qq
U_{z^n,z^{n'},z^{n''}}\ =\ \cases{\,1\qquad {\rm for}\ \ (n,n',n'')\not=(1,1,1)\,,\cr
(-1)^{\sk\,r}\quad\ {\rm for}\ \ n=n'=n''=1\,.}
\qqq
For $\,\sk\,$ and $\,r\,$ odd, the cocycle is cohomologically non-trivial.
As a result
\qq
\sk\ =\ \cases{\,1\qquad{\rm for}\ \ r\ \ {\rm even}\,,\cr
\,2\qquad{\rm for}\ \ r\ \ {\rm odd}\,}
\qqq
and for those values one may take $\,u_{z^{n},z^{n'}}\equiv1\,$ as 
the solution of eq.\,\,(\ref{rtc1}). 
\vskip 0.15cm

\subsection{Groups $\,D_r=Spin(2r)\,$}
\vskip 0.1cm

The Lie algebra of $\,D_r\,$ is $\,so(2r)\,$ composed of imaginary
antisymmetric $\,(2r)\times(2r)$ matrices. The Cartan algebra
may be taken as composed of $\,r\,$ blocks 
$\,\left(\matrix{_0&_{-it_i}\cr^{it_i}&^0}\right)\,$ placed diagonally.
In particular, let $\,e_i\,$ denote the matrix corresponding to $\,t_j=
\delta_{ij}$. \,With the invariant form normalized so that 
$\,\tr\,e_ie_j=\delta_{ij}$, \,roots and coroots of $\,so(2r)\,$ have 
the form $\pm e_i\pm e_j\,$ for $\,i\not=j$. \,The simple roots may be 
chosen as $\,\alpha_i=e_i-e_{i+1}\,$ 
for $\,i=1,\dots,r-1\,$ and $\,\alpha_r=e_{r-1}+e_r$.
%\,The coroot lattice 
%is composed of the elements $\,\sum\limits_{i=1}^rn_ie_i\,$ for integer 
%$\,n_i\,$ such that their sum is even.
\vskip 0.4cm

\noindent I. \ Case of $\,r\,$ odd. 
\vskip 0,3cm

Here the center is $\,\NZ_4\,$ and it may be generated 
by $\,z=e^{-2\pi i\theta}\,$
with $\,\theta=\lambda_r^{\hspace{-0.03cm}^\vee}
={1\over2}\sum\limits_{i=1}^re_i$. \,The permutation $\,z0=r-1,\,z1=r$, 
$\,zi=r-i\,$ for $\,i=2,\dots,r-2$, $\,z(r-1)=1,\,zr=0\,$
induces the extended Dynkin diagram symmetry:

\leavevmode\epsffile[-80 -20 240 165]{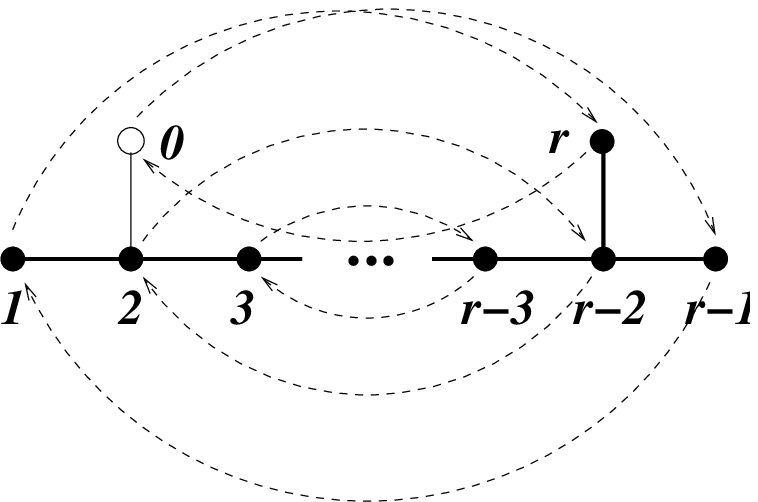}

\noindent $\,SO(2r)=Spin(2r)/\{1,z^2\}$. \,The adjoint action
of $\,w_z\,$ on the Cartan algebra is given by
\qq
w_z\,e_i\,w_z^{-1}\ =\ \cases{\hbox to 1.7cm{$\,\,e_r$\hfill}{\rm for}\ \  
i=1\,,\cr
\hbox to 1.7cm{$-e_{r-i+1}$\hfill}{\rm for}\ \ i\not=1\,.}
\qqq
It may be generated by the product
\qq
r_{\alpha_{r-1}}r_{\alpha_{r-2}}r_{\alpha_{r}}\cdots\cdots\, 
r_{\alpha_4}\cdots\, r_{\alpha_{r-1}}
r_{\alpha_3}\cdots\, r_{\alpha_{r-2}}r_{\alpha_r} 
r_{\alpha_2}\cdots\, r_{\alpha_{r-1}}r_{\alpha_1}\cdots\, r_{\alpha_{r-2}}
r_{\alpha_r}\quad
\qqq
of $\,{r(r-1)\over2}\,$ reflections in simple roots. 
Element $\,w_z\,$ may be taken as a lift to $\,Spin(2r)\,$ of the matrix 
\qq
\left(\matrix{_0&_0&\dots&_0&_1\cr _0&_0&\dots&_1&_0\cr & &\qquad\cdot& & \cr
 & &\ \cdot\ & & \cr & &\cdot\qquad& & \cr _0&_1&{\cdots}&_0&_0\cr
_{-1}&_0&{\cdots}&_0&_0}\right)
\qqq
in $\,SO(2r)$. \,We shall take $\,w_{z^n}=w_z^n\,$ for $\,n=0,1,2,3$. 
\,Then
\qq
c_{z^n,z^m}\ =\ \cases{\hbox to 1.5cm{$\,1$\hfill}{\rm if}\ \ n+m<4\,,\cr
\hbox to 1.5cm{$\,w_z^4$\hfill}{\rm if}\ \ n+m\geq4\,.}
\qqq
It suffices then to determine the value of $\,w_z^4$. \,Since this element
projects to identity in $\,SO(2r)$, \,it is either equal to $\,1\,$ or
to $\,z^2$. \,To determine which is the case, note that we may set 
$\,w_z=\CO\,\ee^{2\pi iX}\CO^{-1}$, \,where $\,\CO\,$ is an element of 
$\,Spin(2r)\,$ projecting to the matrix
\qq
{_1\over^{\sqrt{2}}}\left(\matrix{_{\sqrt{2}}&_0&_0&_0&..&_0&_0&..&_0&_0\cr
_0&_0&_1&_0&..&_0&_0&..&_0&_1\cr_0&_0&_0&_1&..&_0&_0&..&_1&_0\cr & & & &
\ .&.&.&.\ & & &\cr
_0&_0&_0&_0&..&_1&_1&..&_0&_0\cr_0&_0&_0&_0&..&_1&_{-1}&..&_0&_0\cr 
& & & &\ .&.&.&.\ & & &\cr_0&_0&_0&_1&..&_0&_0&..&_{-1}&_0\cr 
_0&_0&_1&_0&..&_0&_0&..&_0&_{-1}\cr_0&_{\sqrt{2}}&_0&_0&..&_0&_0&..&_0&_0}
\hspace{-0.4cm}\right)
\qqq
and $\,\,X={1\over4}e_1+{1\over2}(e_{r+3\over2}\dots+e_r)\,\,$
so that $\,\ee^{2\pi iX}\,$ projects to the matrix
\qq
\hspace{-1cm}
\ee^{2\pi iX}\ =\ \left(\matrix{_0&_1&_0&_0&..&_0&_0&..&_0&_0\cr
_{-1}&_0&_0&_0&..&_0&_0&..&_0&_0\cr_0&_0&_1&_0&..&_0&_0&..&_0&_0\cr
_0&_0&_0&_1&..&_0&_0&..&_0&_0\cr & & & &..& & &..& & \cr
_0&_0&_0&_0&..&_1&_0&..&_0&_0\cr_0&_0&_0&_0&..&_0&_{-1}&..&_0&_0\cr
 & & & &..& & &..& & \cr_0&_0&_0&_0&..&_0&_0&..&_{-1}&_0\cr
_0&_0&_0&_0&..&_0&_0&..&_0&_{-1}}\right)
\qqq
in $\,SO(2r)$. \,It follows that $\,w_z^4=e^{\m8\pi i X}=z^2\,$ since 
$\,4X\,$ is not in the coroot lattice. We may take 
\qq
e_{z^n,z^m}\ =\ \cases{\hbox to 1.5cm{$\,\,\,0$\hfill}{\rm if}\ \ n+m<4\,,\cr
\hbox to 1.5cm{$\,2\theta$\hfill}{\rm if}\ \ n+m\geq4\,.}
\qqq
This results in
\qq
U_{z^n,z^{n'},z^{n''}}\ =\ (-1)^{\,\sk\,n''{n+n'-[n+n']\over 4}}
\label{cdr}
\qqq
for $\,0\leq n,n',n''\leq3$, \,where now $\,[m]=m\,\,{\rm mod}\,\,4\,$
with $\,0\leq[m]\leq3$. $\,U\,$ is cohomologically non-trivial for
$\,\sk\,$ odd, hence $\,\sk=2\,$ if $\,Z=\NZ_4$. \,On the other hand,
the cocycle (\ref{cdr}) becomes trivial when restricted to the cyclic
subgroup of order 2 generated by $\,z^2\,$ so that $\,\sk=1\,$ if
$\,Z=\NZ_2$. \,In both cases, for the above values of $\,\sk$, \,one 
may take $\,u_{z^{n},z^{n'}}\equiv1\,$ as the solution of 
eq.\,\,(\ref{rtc1}). 
\vskip 0.3cm

\noindent II. \,\,Case of $\,r\,$ even. 
\vskip 0,3cm

Here the center is $\,\NZ_2\times\NZ_2$. It is generated by $\,z_1=
\ee^{-2\pi i\theta_1}\,$ and 
$\,z_2=\ee^{-2\pi i\theta_2}\,$ for 
$\,\theta_1=\lambda_r^{\hspace{-0.03cm}^\vee}
={1\over2}(\sum\limits_{i=1}^re_i)\,$
and $\,\theta_2=\lambda_1^{\hspace{-0.03cm}^\vee}=e_1$.
\,These elements induce the permutations 
$\,z_10=r,\,z_1i=r-i\,$ for $\,i=1,\dots,r-1$,
$\,z_1r=0$, $\,z_20=1,\,z_21=0,\,z_2i=i\,$ for $\,i=2,\dots,r-2$,
$\,z_2(r-1)=r,\,z_2r=r-1\,$ giving rise to the symmetries
of the extended Dynkin diagrams:

\leavevmode\epsffile[-80 -10 240 160]{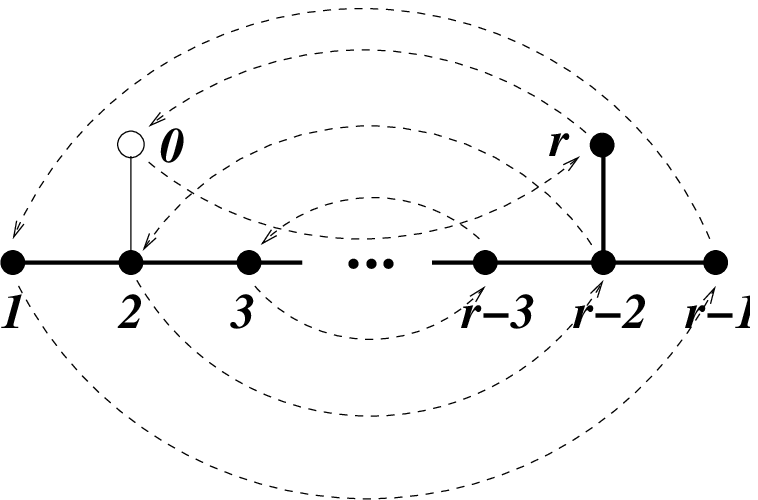}

\leavevmode\epsffile[-80 -20 240 80]{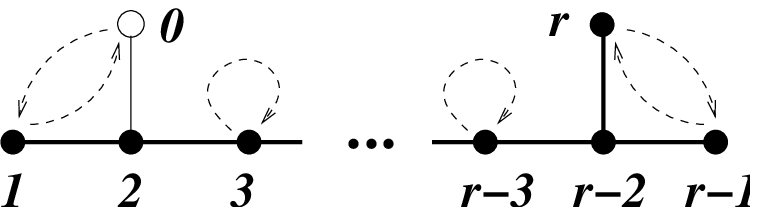}

\noindent$\,SO(2r)=Spin(2r)/\{1,z_2\}$. \,The adjoint actions of 
$\,w_{z_1}\,$ and $\,w_{z_2}\,$ on the Cartan algebra are given by
\qq
w_{z_1}\,e_i\,w_{z_1}^{-1}\ =\ -e_{r-i+1}\,.\qquad
w_{z_2}\,e_i\,w_{z_2}^{-1}\ =\ \cases{\hbox to 1.5cm{$-e_i$\hfill}{\rm for}
\ \ i=1,r\,,\cr \hbox to 1.5cm{$\,\,e_i$\hfill}{\rm for}\ \ i\not=1,r\,.}
\qqq
They may be generated by the products
\qq
r_{\alpha_{r}}\cdots\cdots\, 
r_{\alpha_4}\cdots\, r_{\alpha_{r-1}}
r_{\alpha_3}\cdots\, r_{\alpha_{r-2}}r_{\alpha_r} 
r_{\alpha_2}\cdots\, r_{\alpha_{r-1}}r_{\alpha_1}\cdots\, r_{\alpha_{r-2}}
r_{\alpha_r}\ 
\qqq
and
\qq
r_{\alpha_1}\cdots r_{\alpha_{r-2}}r_{\alpha_r}
r_{\alpha_{r-1}}\cdots r_{\alpha_2}r_{\alpha_1}
\qqq
of, respectively, $\,{r(r-1)\over2}\,$ and $\,2(r-1)\,$ reflections 
in simple roots. 
Elements $\,w_{z_1}\,$ and $\,w_{z_2}\,$ may be taken as lifts to 
$\,Spin(2r)\,$ of the $\,SO(2r)\,$ matrices 
\qq
\left(\matrix{_0&_0&_0&\dots&_0&_0&_1\cr_0&_0&_0&\dots&_0&_1&_0\cr
_0&_0&_0&\dots&_1&_0&_0\cr.&.&.&\dots&.&.&.\cr.&.&.&\dots&.&.&.
\cr_0&_0&_1&\dots&_0&_0&_0\cr_0&_1&_0&\dots&_0&_0&_0\cr
_1&_0&_0&\dots&_0&_0&_0}\right)\quad\ {\rm and}\ \quad
\left(\matrix{_{-1}&_0&_0&\dots&_0&_0&_0\cr_0&_1&_0&\dots&_0&_0&_0\cr
_0&_0&_1&\dots&_0&_0&_0\cr.&.&.&\dots&.&.&.\cr
.&.&.&\dots&.&.&.\cr_0&_0&_0&\dots&_1&_0&_0\cr_0&_0&_0&\dots&_0&_1&_0\cr
_0&_0&_0&\dots&_0&_0&_{-1}}\right)
\qqq
respectively. We may set 
\qq
w_{z_1}\ =\ \CO_1\,\ee^{2\pi i X_1}\,\CO_1^{-1}\,,\qquad
w_{z_2}\ =\ \CO_2\,\ee^{2\pi i X_2}\,\CO_2^{-1}\ =\ 
\CO_1\CO_2\,\ee^{2\pi i X_2}\,\CO_2^{-1}\CO_1^{-1}\quad
\qqq
where $\,\CO_i\,$ are in $\,Spin(2r)\,$ and project to the
$\,SO(2r)\,$ matrices
\qq
\hspace{-0.8cm}{_1\over^{\sqrt{2}}}\left(\matrix{_1&_0&..&_0&_0&..&_0&_1\cr
_0&_1&..&_0&_0&..&_1&_0\cr & &
\ .&.&.&.\ & & &\cr
_0&_0&..&_1&_1&..&_0&_0\cr_0&_0&..&_1&_{-1}&..&_0&_0\cr 
& &\ .&.&.&.\ & & &\cr_0&_1&..&_0&_0&..&_{-1}&_0\cr 
_1&_0&..&_0&_0&..&_0&_{-1}}
\hspace{-0.4cm}\right)\ \quad{\rm and}\quad\ 
\left(\matrix{_1&_0&_0&_0&..&_0&_0\cr_0&_0&_1&_0&..&_0&_0\cr
_0&_0&_0&_1&..&_0&_0\cr.&.&.&.&..&.&.\cr.&.&.&.&..&.&.\cr
_0&_0&_0&_0&..&_1&_0\cr_0&_0&_0&_0&..&_0&_1\cr_0&_1&_0&_0&..&_0&_0}
\right),
\qqq
respectively, with
$\,\,X_1={1\over2}(e_{{r\over 2}+1}+\dots+e_r)\,\,$ and $\,\,
X_2={1\over2}e_1$. \,The exponentials $\,\ee^{2\pi i X_1}\,$ 
and $\,\ee^{2\pi i X_2}\,$ project in turn to the matrices
\qq
\left(\matrix{_1&_0&..&_0&_0&..&_0&_0\cr_0&_1&..&_0&_0&..&_0&_0\cr
.&.&..&.&.&..&.&.\cr_0&_0&..&_1&_0&..&_0&_0\cr
_0&_0&..&_0&_{-1}&..&_0&_0\cr.&.&..&.&.&..&.&.\cr
_0&_0&..&_0&_0&..&_{-1}&_0\cr_0&_0&..&_0&_0&..&_0&_{-1}}\right)
\ \quad{\rm and}\quad\ \left(\matrix{_{-1}&_0&_0&_0&..&_0&_0\cr
_0&_{-1}&_0&_0&..&_0&_0\cr_0&_0&_1&_0&..&_0&_0\cr
_0&_0&_0&_1&..&_0&_0\cr.&.&.&.&..&.&.\cr.&.&.&.&..&.&.\cr
_0&_0&_0&_0&..&_1&_0\cr_0&_0&_0&_0&..&_0&_1}\right)
\qqq
respectively. Since $\,w_{z_i}^2\,$ projects to $\,1\,$ in $\,SO(2r)\,$ it
is equal to $\,1\,$ or to $\,z_2\,$ in $\,Spin(2r)$. \,Which is the case,
depends on whether $\,2X_i\,$ is in the coroot lattice. We infer that
\qq
&&w_{z_1}^2\ =\ \cases{\hbox to 1cm{$\,1$\hfill}{\rm if}\ \ r\ \ {\rm is\ divisible\ by}
\ \,4\,,\cr
\hbox to 1cm{$\,z_2$\hfill}{\rm otherwise},}\cr\cr
&&w_{z_2}^2\ =\ z_2\,.
\qqq
Besides,
\qq
&&w_{z_1}w_{z_2}w_{z_1}^{-1}w_{z_2}^{-1}\ =\ \CO_1\Big(\ee^{2\pi iX_1}
\CO_2\,\ee^{2\pi iX_2}\CO_2^{-1}\ee^{-2\pi iX_1}\CO_2\,
\ee^{-2\pi iX_2}\CO_2^{-1}\Big)\CO_1^{-1}\cr
&&=\CO_1\Big(\ee^{2\pi iX_1}w_{z_2}\,\ee^{-2\pi iX_1}w_{z_2}^{-1}\Big)
\CO_1^{-1}\ =\ \CO_1\,\ee^{2\pi i e_r}\CO_1^{-1}\ =\ z_2\,.
\qqq
Setting $\,w_1=1\,$ and $\,w_{z_1z_2}=w_{z_1}w_{z_2}$, \,we infer that
for $\,r\,$ divisible by $\,4$,
\qq
c_{z,z'}\,=\,\cases{\,\hbox to 0.9 cm{$z_2$\hfill}{\rm if}\ \ (z,z')=(z_2,z_1),\,
(z_2,z_2),\,(z_1z_2,z_1),\,(z_1z_2,z_2)\,,\cr\cr
\,\hbox to 0.9 cm{$1$\hfill}{\rm otherwise}}
\qqq
and for $\,r\,$ not divisible by $\,4$,
\qq
c_{z,z'}\,=\,\cases{\,\hbox to 0.9 cm{$z_2$\hfill}{\rm if}\ \ 
(z,z')=(z_1,z_1),\,
(z_1,z_1z_2),\,(z_2,z_1),\cr
\hspace{3cm}(z_2,z_2),\,(z_1z_2,z_2),\,(z_1z_2,z_1z_2)\,,\cr\cr
\,\hbox to 0.9 cm{$1$\hfill}{\rm otherwise}\,.}\qquad\ \ 
\qqq
We may then take for $\,r\,$ divisible by $\,4$,
\qq
e_{z,z'}\,=\,\cases{\,\hbox to 0.9 cm{$\theta_2$\hfill}{\rm if}\ \ (z,z')=(z_2,z_1),\,
(z_2,z_2),\,(z_1z_2,z_1),\,(z_1z_2,z_2)\,,\cr\cr
\,\hbox to 0.9 cm{$0$\hfill}{\rm otherwise}\,,}
\qqq
and for $\,r\,$ not divisible by $\,4$,
\qq
e_{z,z'}\,=\,\cases{\,\hbox to 0.9 cm{$\theta_2$\hfill}{\rm if}\ \ 
(z,z')=(z_1,z_1),\,
(z_1,z_1z_2),\,(z_2,z_1),\cr
\hspace{3cm}(z_2,z_2),\,(z_1z_2,z_2),\,(z_1z_2,z_1z_2)\,,\cr\cr
\,\hbox to 0.9 cm{$0$\hfill}{\rm otherwise}\,.}\qquad\ \ 
\qqq
Explicit calculation gives:
\qq
U_{z,z',z''}\ =\ \cases{\,(-1)^{\sk(1+r/2)}\quad\,{\rm for}\ \ (z,z',z'')=(z_1z_2,z_1,z_1),\,(z_1z_2,z_1,z_1z_2)\,,\cr
\,(-1)^{\sk\m r/2}\quad\ \ \quad{\rm for}\ \ (z,z',z'')=(z_1,z_1,z_1),\,(z_1,z_1,z_1z_2),\,(z_1,z_1z_2,z_1),\cr
\hspace{2.5cm}\,(z_1,z_1z_2,z_1z_2),\,(z_1z_2,z_1z_2,z_1),
\,(z_1z_2,z_1z_2,z_1z_2)\,,\cr
\,(-1)^\sk\,\m\ \ \qquad\quad{\rm for}\ \ (z,z',z'')=(z_2,z_1,z_1),
\,(z_2,z_1,z_1z_2),\,(z_2,z_2,z_1),\cr
\hspace{2.5cm}\,(z_2,z_2,z_1z_2),\,(z_1z_2,z_2,z_1),\,(z_1z_2,z_2,z_1z_2)
\,,\cr
\,1\hspace{2.3cm}{\rm otherwise}.}
\qqq
The cocycle $\,U\,$ is cohomologically nontrivial if $\,\sk\m r/2\,$
is odd. If $\,\sk\,$ is even, it is trivial, and any 2-cocycle $\,u\,$ solves
eq.\,\,(\ref{rtc1}). In particular, we may take
\qq
u_{z,z'}\ =\ \cases{\,\hbox to 1cm{$\pm1$\hfill}{\rm for}\ \ 
(z,z')=(z_2,z_1),\,(z_2,z_1z_2),\,(z_1z_2,z_1),\,(z_1z_2,z_1z_2)\,,\cr
\,\hbox to 1cm{$1$\hfill}{\rm otherwise}}
\label{necc}
\qqq
representing two non-equivalent classes in $\,H^2(Z,U(1))$.
\,When $\,\sk\,$ is odd and $\,r/2\,$ is even then $\,U\,$ is 
cohomologically trivial and 
\qq
u_{z,z'}\ =\ \cases{\,\hbox to 1cm{$\pm i$\hfill}{\rm for}\ \ 
(z,z')=(z_2,z_1),\,(z_2,z_1z_2),\,(z_1z_2,z_1),\,(z_1z_2,z_1z_2)\,,\cr
\,\hbox to 1cm{$1$\hfill}{\rm otherwise}}
\label{nech}
\qqq
give two solutions of eq.\,\,(\ref{rtc1}) differing by a nontrivial cocycle
(\ref{necc}). Hence for $\,Z=\NZ_2\times\NZ_2$, $\,\sk=1\,$ if $\,r/2\,$ is even
and $\,\sk=2\,$ for $\,r/2\,$ odd. 
\vskip 0.3cm

If $\,Z\,$ is the $\NZ_2\,$ subgroup generated by $\,z_1\,$
or by $\,z_1z_2\,$ then the restriction of $\,U\,$ to $\,Z\,$ is
cohomologically nontrivial if $\,\sk\, r/2\,$ is odd and is trivial
if $\,\sk\m r/2\,$ is even. Hence $\,\sk=1\,$ if $\,r/2\,$ is even
and $\,\sk=2\,$ if it is odd. For $\,Z=\NZ_2\,$ generated by $\,z_2$,
\,the restriction of $\,U\,$ to $\,Z\,$ is trivial so that $\,\sk=1$.
One may take $\,u_{z,z'}\equiv1\,$ as the solution of 
eq.\,\,(\ref{rtc1}) in those cases. 
\vskip 0.2cm

\subsection{Group $\,E_6\,$}
\vskip 0.15cm

We shall identify the Cartan algebra of the exceptional group $\,E_6\,$
with the subspace of $\,\NR^7\,$ with the first six coordinates
summing to zero, with the scalar product inherited from $\,\NR^7$.
\,The simple roots, may be taken as $\,\alpha_i=e_i-e_{i+1}\,$ for
$\,i=1,\dots,5\,$ and $\,\alpha_6={1\over2}(-e_1-e_2-e_3+e_4+e_5+e_6)+
{1\over\sqrt{2}}e_7$, \,where $\,e_i\,$ are the vectors of the canonical
bases of $\,\NR^7$. \,The center of $\,E_6\,$ is $\,\NZ_3\,$ and
it is generated by $\,z=\ee^{-2\pi i\theta}\,$ with $\,\theta
=\lambda_5^{\hspace{-0.03cm}^\vee}={1\over6}
(e_1+e_2+e_3+e_4+e_5-5e_6)+{1\over\sqrt{2}}e_7$. 
\,The permutation $\,z0=1,\,z1=5,z2=4,\,z3=3,\,z4=6,
\,z5=0,\,z6=2\,$ induces the symmetry of the extended Dynkin diagram:

\leavevmode\epsffile[-100 -20 220 170]{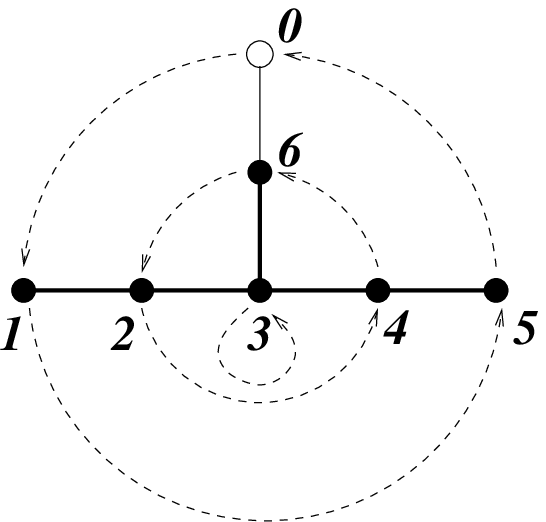}

\noindent The adjoint action
of $\,w_z\,$ on the Cartan algebra may be generated by setting
\qq
&&w_z\,e_1\,w_z^{-1}\ =\ -e_6\,,\qquad w_z\,e_2\,w_z^{-1}\ =\ -e_5\,,\qquad 
w_z\,e_3\,w_z^{-1}\ =\ -e_4\,,\cr
&&w_z\,e_4\,w_z^{-1}\ =\ -e_3\,,\qquad 
w_z\,e_5\,w_z^{-1}\ =\ {_1\over^2}(e_1+e_2-e_3-e_4-e_5-e_6)-{_1\over^{\sqrt{2}}}e_7\,,\cr
&&w_z\,e_6\,w_z^{-1}\ =\ {_1\over^2}(e_1+e_2-e_3-e_4-e_5-e_6)+{_1\over^{\sqrt{2}}}e_7\,,\cr
&&w_z\,e_7\,w_z^{-1}\ =\ {_1\over^{\sqrt{2}}}(-e_1+e_2)
\qqq
and is given by the product
\qq
r_{\alpha_1}r_{\alpha_2}r_{\alpha_3}r_{\alpha_4}r_{\alpha_5}r_{\alpha_6}
r_{\alpha_3}r_{\alpha_2}r_{\alpha_1}r_{\alpha_4}r_{\alpha_3}r_{\alpha_2}
r_{\alpha_6}r_{\alpha_3}r_{\alpha_4}r_{\alpha_5}\label{prref}
\qqq
of $\,16\,$ reflections that may be rewritten as the product of
4 reflections $\,r_{\beta_1}r_{\beta_4}r_{\beta_5}r_{\beta_2}\,$
in non-simple roots
\qq
&&\beta_1\ =\ \alpha_1+\alpha_2+\alpha_3+\alpha_4\,,\qquad
\beta_2\ =\ \alpha_3+\alpha_4+\alpha_5+\alpha_6\,,\cr
&&\beta_4\ =\ \alpha_1+\alpha_2+\alpha_3+\alpha_6\,,\qquad
\beta_5\ =\ \alpha_2+\alpha_3+\alpha_4+\alpha_5\,,
\qqq
The family of roots $\,(\beta_1,\beta_2,\beta_3,\beta_4,\beta_5,\beta_6)\,$
with
\qq
&&\beta_3\ =\ -\alpha_1-\alpha_2-2\alpha_3-\alpha_4-\alpha_5-\alpha_6\,,
\qquad\beta_6\ =\ \alpha_3
\qqq
provides another set of simple roots for $\,E_6\,$ corresponding 
to the same Cartan matrix. The roots $\,\beta_i\,$ with $\,i\leq 5\,$ 
and their step generators $\,e_{\pm\beta_i}\,$ generate an $\,A_5\,$ 
subalgebra of $\,E_6\,$ which, upon exponentiation, gives rise to 
an $\,SU(6)\,$ subgroup of group $\,E_6$. \,The group elements that 
implement by conjugation the Weyl reflections $\,r_{\beta_i}\,$ 
of the Cartan algebra of $\,E_6\,$ may be taken as 
$\,\ee^{\,{\pi\over2i}(e_{\beta_i}+e_{-\beta_i})}\,$ so that
they belong to the $\,SU(6)\,$ subgroup for $\,i\leq 5$.
\,We infer that, identifying roots $\,\beta_i\,$ for $\,i\leq 5\,$
with the standard roots of $\,A_5$, \,the element $\,w_z\,$ may be
taken as the matrix
\qq
\left(\matrix{_0&_0&_1&_0&_0&_0\cr _1&_0&_0&_0&_0&_0\cr
_0&_1&_0&_0&_0&_0\cr _0&_0&_0&_0&_0&_1\cr _0&_0&_0&_1&_0&_0\cr
_0&_0&_0&_0&_1&_0}\right)\ \ \in\ SU(6)\,\subset\,E_6
\qqq
which satisfies $\,w_z^3=1$. Setting now $\,w_1=1\,$ and $\,w_{z^2}=w_z^2$, we 
end up with trivial cocycle $\,c_{z,z'}\,$ corresponding to 
$\,e_{z,z'}\equiv0$.
\,Consequently, the cocycle $\,U\,$ is also trivial, $\,\sk=1\,$ and 
$\,u_{z^n,z^{n'}}\equiv1\,$ solves eq.\,\,(\ref{rtc1}).
\vskip 0.2cm

\subsection{Group $\,E_7\,$}
\vskip 0.1cm

The Cartan algebra of $\,E_7\,$ may be identified with the subspace
of $\,\NR^8\,$ orthogonal to the vector $\,(1,1,\dots,1)\,$ with
the simple roots $\,\alpha_i=e_i-e_{i+1}\,$ for $\,i=1,\dots,6\,$
and $\,\alpha_7={1\over2}(-e_1-e_2-e_3-e_4+e_5+e_6+e_7+e_8)\,$
with $\,e_i\,$ the vectors of the canonical basis of $\,\NR^8$.
\,The center of $\,E_7\,$ is $\,\NZ_2\,$ with the non-unit element
$\,z=\ee^{-2\pi i\,\theta}\,$ for $\,\theta
=\lambda_1^{\hspace{-0.03cm}^\vee}={1\over 4}(3,-1,-1,-1,-1,-1,-1,3)$. 
\,The permutation $\,z0=1,\,z1=0,\,z2=6,\,z3=5,\,z4=4,\,z5=3,\,z6=2,\,z7=7\,$
generates the symmetry of the extended Dynkin diagram:

\leavevmode\epsffile[-76 -20 244 170]{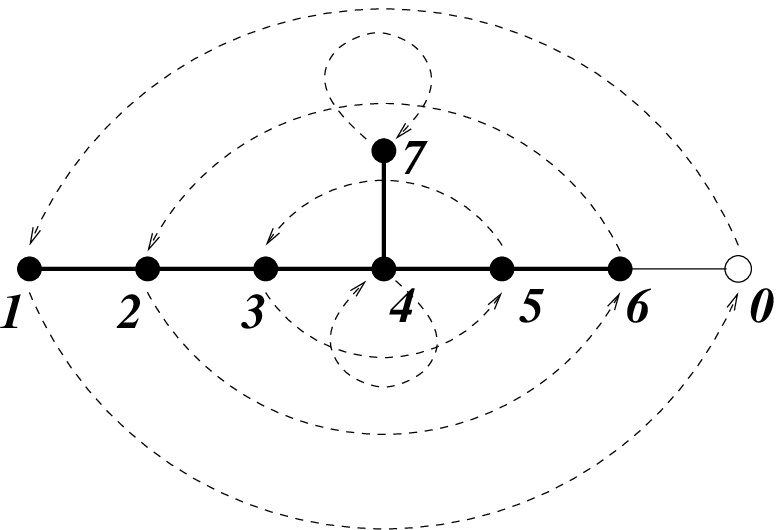}

\noindent The adjoint action of $\,w_z\,$ may by obtained by setting
\qq
w_z\,e_i\,w_z^{-1}\ =\ -e_{9-i}
\qqq
and is given by the product
\qq
&r_{\alpha_1}r_{\alpha_2}r_{\alpha_3}r_{\alpha_4}r_{\alpha_5}r_{\alpha_7}
r_{\alpha_4}r_{\alpha_6}r_{\alpha_3}r_{\alpha_5}r_{\alpha_2}r_{\alpha_4}
r_{\alpha_1}r_{\alpha_3}r_{\alpha_7}&\cr
&\cdot\,r_{\alpha_4}r_{\alpha_2}r_{\alpha_5}
r_{\alpha_3}r_{\alpha_6}r_{\alpha_4}r_{\alpha_7}r_{\alpha_5}r_{\alpha_4}
r_{\alpha_3}r_{\alpha_2}r_{\alpha_1}&\label{prref7}
\qqq
of $\,27\,$ simple root reflections that may be rewritten as the product of
3 reflections $\,r_{\beta_1}r_{\beta_3}r_{\beta_7}\,$ for 
\qq
&&\beta_1\ =\ \alpha_1+2\alpha_2+2\alpha_3+2\alpha_4+\alpha_5+\alpha_7\ 
=\ \omega(\alpha_1)\,,\cr
&&\beta_3\ =\ \alpha_1+\alpha_2+2\alpha_3+2\alpha_4+\alpha_5+\alpha_6
+\alpha_7\ =\ \omega(\alpha_3)\,,\cr
&&\beta_7\ =\ \alpha_1+\alpha_2+\alpha_3+2\alpha_4+2\alpha_5+
\alpha_6+\alpha_7\ =\ \omega(\alpha_7)\,,
\qqq
where $\ \omega=r_{\alpha_1}r_{\alpha_2}r_{\alpha_3}r_{\alpha_4}r_{\alpha_5}r_{\alpha_7}
r_{\alpha_4}r_{\alpha_6}r_{\alpha_3}r_{\alpha_5}r_{\alpha_2}r_{\alpha_4}$.
\,The roots $\,\beta_1,\,\beta_3\,$ and $\,\beta_7\,$ may be completed to a new system 
of simple roots of $\,E_7\,$ by setting
\qq
&&\beta_2\ =\ -(\alpha_1+\alpha_2+2\alpha_3+2\alpha_4+\alpha_5+\alpha_7)\ 
=\ \omega(\alpha_2)\,,\cr
&&\beta_4\ =\ -(\alpha_1+\alpha_2+\alpha_3+2\alpha_4+\alpha_5+\alpha_6
+\alpha_7)\ =\ \omega(\alpha_4)\,,\cr
&&\beta_6\ =\ \alpha_7\ =\ \omega(\alpha_6)\,.
\qqq
In particular, $\,\beta_1,\beta_2,\beta_3,\beta_4,\beta_7\,$ and their step generators
span a subalgebra $\,A_5\subset E_7\,$ that, upon exponentiation, gives rise to a subgroup
$\,SU(6)\,$ in group $\,E_7$. \,The element $\,w_z\,$ implementing by conjugation
the Weyl transformation (\ref{prref7}) may be chosen as
\qq
\left(\matrix{_0&_1&_0&_0&_0&_0\cr _{-1}&_0&_0&_0&_0&_0\cr
_0&_0&_0&_1&_0&_0\cr _0&_0&_{-1}&_0&_0&_0\cr _0&_0&_0&_0&_0&_1\cr
_0&_0&_0&_0&_{-1}&_0}\right)\ \ \in\ SU(6)\,\subset\,E_7
\qqq
upon identifying of the roots $\,\beta_1,\beta_2,\beta_3,\beta_4,\beta_7\,$ with the 
standard roots of $\,A_5=su(6)$. \,In particular, $\,w_z^2=-1\in SU(6)\,$ or
\qq
w_z^2\ =\ \ee^{\,\pi i(\beta_1+\beta_3+\beta_7)}\ =\ (-i,-i,-i,-i,-i,-i,-i,-i)\ =\ 
\ee^{\,2\pi i\,\theta}\,.
\qqq
With that choice of $\,w_z$, we infer that
\qq
c_{1,1}\ =\ c_{1,z}\ =\ c_{z,1}\ =\ 1\,,\qquad c_{z,z}\ =\ w_z^2\ =\ 
\ee^{\,2\pi i\,\theta}
\quad
\qqq
and we may take
\qq
e_{1,1}\ =\ e_{1,z}\ =\ e_{z,1}\ =\ 0\,,\qquad e_{z,z}\ =\ \theta\,.
\qqq
This leads to
\qq
U_{z^n,z^{n'},z^{n''}}\ =\ \cases{\,1\qquad\ {\rm for}\ \ (n,n',n'')\not=(1,1,1)\,,\cr
(-1)^{\sk}\qquad\ {\rm for}\ \ n=n'=n''=1\,.}
\qqq
$U\,$ is trivial if $\,\sk\,$ is even and is cohomologically nontrivial when
$\,\sk\,$ is odd. Hence $\,\sk=2\,$ and one may take
$\,u_{z^{n},z^{n'}}\equiv1\,$ as the solution of eq.\,\,(\ref{rtc1})
for that value of $\,\sk$. 
\vskip 0.3cm

\nsection{Conclusions}

We have presented an explicit construction of the basic
gerbes over groups $\,G'=G/Z\,$ where $\,G\,$ is a simple compact 
connected and simply connected group and $\,Z\,$ is a 
subgroup of the center of $\,G$. \,By definition of the basic gerbe, the 
pullback to $\,G\,$ of its curvature $\,H'\,$ is the closed 3-form $\,H=
{\sk\over{12\pi}}\,\tr\,(g^{-1}dg)^3\,$ with the level $\,\sk\,$
taking the lowest possible positive value. The restriction on
$\,\sk\,$ came from the cohomological equation (\ref{rtc1}) 
that assures the associativity of the gerbe's groupoid product.
In agreement with the general theory, see \cite{Gaw0,GR}, the levels $\,\sk\,$ 
of the basic gerbes are the lowest positive numbers for which the periods 
of $\,H'\,$ belong to $\,2\pi\NZ$. \,They have been previously found in 
ref.\,\,\cite{FGK} and we have recovered here the same set of numbers. 
The basic gerbe over $\,G'\,$ is unique up to stable isomorphisms except 
for $\,G'=Spin(4N)/\NZ_2\times \NZ_2$. In the latter case,
using the two different choices of sign in the solutions (\ref{necc}) 
or (\ref{nech}) of the cohomological relation (\ref{rtc1}), one obtains 
basic gerbes belonging to two different stable isomorphism classes,
the doubling already observed in ref.\,\,\cite{FGK}. We plan to use 
the results of the present paper in order to extend the classification 
of the fully symmetric branes in groups $\,SU(N)/Z\,$ worked out 
in ref.\,\,\cite{GR} to all groups $\,G'$.
\vskip -0.2cm

\nappendix{A}
\vskip 0.4cm

\noindent We shall obtain here the condition (\ref{rtc}) for
the associativity of the groupoid product $\,\mu'\,$ defined
by (\ref{map}). Let $\,(y,y',y'',y''')\in{Y'}^{[4]}\,$ with
$\,(y,zy',z(z'y''),z(z'(z''y''')))\,$ belonging to
$\,Y_{ijkl}\,$ and projecting to $\,g\in\CO_{ijkl}$.
Taking $\,h\in G\,$ such that $\,\rho_{ijkl}(g)=hG_{ijkl}$,
we may complete eqs.\,\,(\ref{wmr0}) to (\ref{wmr2}) by
\qq
&&z(z'(z''y'''))\ =\ (g\m,\s h\,{{\gamma}'''}^{-1})\,,\label{wmr3}\\ \cr
&&y'''\ =\ ((zz'z'')^{-1}g\m,\, h\,w_zw_{z'}w_{z''}
\,((\gamma_z''')_{z'})_{z''}^{-1})\,,\label{wmr4}\\ \cr
&&z''y'''\ =\ ((zz')^{-1}g\m,\, h\,w_zw_{z'}\,
(\gamma_z''')_{z'}^{-1})\,,\label{wm5}\\ \cr
&&(z'z'')y'''\ =\ (z^{-1}g\m,\, h\,w_z\,(c_{z',z''}^{\,-1}\gamma_z''')^{-1})\,,
\label{wm6}\\ \cr
&&(zz'z'')y'''\ =\ (g\m,\,h\,(c_{zz',z''}^{\,-1}c_{z,z'}^{\,-1}
{\gamma'''})^{-1})\,.
\label{wm7}
\qqq
The $\,G_{ijk}$-orbits in eqs.\,\,(\ref{ell1}) to (\ref{ell3}) may be now
replaced by the $\,G_{ijkl}$-orbits.
We shall need further line-bundle elements. Let
\qq
&&\ell_{k_{zz'}l_{zz'}}\ =\,\ ((zz')^{-1}g\m,\,h\, w_zw_{z'}\m,\,
[({\tilde\gamma}''_z)_{z'},
({\tilde\gamma}'''_z)_{z'},\m u'']^{^{\,\sk}}_{_{k_{zz'}l_{zz'}}})
G_{i_{zz'}j_{zz'}k_{zz'}l_{zz'}}\cr\cr
&&=\ \chi_{_{kl}}
(\tilde c_{z,z'})^{^{\sk}}\,((zz')^{-1}g\m,\,h\, w_{zz'}\m,\,
[(\tilde c_{z,z,'}^{\,-1}{\tilde\gamma}'')_{zz'},
(\tilde c_{z,z'}^{\,-1}{\tilde\gamma}''')_{zz'}\m,
\,u'']^{^{\,\sk}}_{_{k_{zz'}l_{zz'}}})
G_{i_{zz'}j_{zz'}k_{zz'}l_{zz'}}\cr\cr
&&\hspace{7cm}\in\,\ 
L^{^{\,\sk}}_{(y'',z''y''')}\,=\,L'_{(y'',y''')}\,,\label{ell4}
\qqq
where we have used the identifications entering the definition of the line
bundle $\,L_{k_{zz'}\ell_{zz'}}$. \,Similarly, let
\qq
&&\ell_{i\m l}\ \ \ \ \ \ =
\,\ (g\m,\, h,\,[{\tilde\gamma},
\,\tilde c_{zz',z''}^{\,-1}\tilde c_{z,z'}^{\,-1}{\tilde\gamma}'''
,\,u\m u'u'' ]^{^{\,\sk}}_{il})G_{ijkl}\cr\cr
&&=\,\ \chi_{_l}((\delta\tilde c)_{z,z',z''})^{^{\sk}}\,(g\m,\, h,\,[{\tilde
\gamma},\,\tilde c_{z,z'z''}^{\,-1}\,w_z(\tilde c_{z',z''}^{\,-1})
w_z^{-1}\,{\tilde
\gamma}''',\,u\m u'u'' ]^{^{\,\sk}}_{il})G_{ijkl}\cr\cr
&&\hspace{7cm}\in\,\ L^{^{\,\sk}}_{(y,(zz'z'')y''')}\,\ =\ 
L'_{(y,y''')}\,,\ \ \ \ \,\qquad\label{ell5}
\qqq
where the 3-cocycle $\,\delta\tilde c\,$ is given by (\ref{wib}). 
Finally, let
\qq
&&\ell_{j_zl_z}\ \ \ \,=\,\ (z^{-1}g\m,\,h\,w_z,\,[{\tilde\gamma}'_z,
\tilde c_{z',z''}^{\,-1}{\tilde\gamma}_z''',\m u'u'']^{^{\,\sk}}_{j_zl_z})
G_{i_zj_zk_zl_z}\cr\cr
&&\hspace{7cm}\in\,\ L^{^{\,\sk}}_{(y',(z'z'')y''')}\,
\,\,=\ L'_{(y',y''')}\,.\ \quad\qquad
\label{ell6}
\qqq
Now
\qq
&\displaystyle{\mu'\left(\mu'\m(\m\ell_{ij}\otimes\m\ell_{j_zk_z})\m\otimes
\ell_{k_{zz'}l_{zz'}}\right)\ =\ u^{ijk}_{z,z'}\,\,\mu'\left(\ell_{ik}
\otimes\m\ell_{k_{zz'}l_{zz'}}\right)}&\cr\cr
&\displaystyle{=\ 
u^{ijk}_{z,z'}\,\,u^{ikl}_{zz',z''}\,\,\chi_{_{kl}}(\tilde c_{z,z'})^{^{\sk}}
\ \ell_{il}\,.}&
\qqq
On the other hand,
\qq
&\displaystyle{\mu'\left(\ell_{ij}\otimes\m\mu'\m(\m\ell_{j_zk_z}\otimes
\ell_{k_{zz'}l_{zz'}})\right)\ =\ u^{j_zk_zl_z}_{z',z''}\,\,\mu'\left(\ell_{ij}
\otimes\m\ell_{j_{z}l_{z}}\right)}&\cr\cr
&\displaystyle{=\ 
u^{j_zk_zl_z}_{z',z''}\,\,u^{ijl}_{z,z'z''}\,\,
\chi_{l}((\delta\tilde c)_{z,z',z''})^{^{-\sk}}
\ \ell_{il}\,.}&
\qqq
Equating both expressions, we infer condition (\ref{rtc}).
\medskip

\hspace{10cm}$\Box$
\vskip -0.46cm

\nappendix{B}
\vskip 0.4cm

\noindent{\bf Proof of Lemma 3.}\ \ With $\,(u^{ijk}_{z,z'})\,$ 
given by eq.\,\,(\ref{tl2}) and
$\,(u_{z,z'})\,$ solving eq.\,\,(\ref{rtc1}), the left hand side
of (\ref{rtc}) becomes
\qq
&&\,\,\,\chi_{_{l_z\m(z'z''z''0)}}(\tilde c_{z',z''})^{^{-\sk}}\,
\chi_{_{l\m(zz'z''0)}}(\tilde c_{zz',z''})^{^{\sk}}\,
\chi_{_{l\m(zz'z''0}}(\tilde c_{z,z'z''})^{^{-\sk}}\cr\cr
&&\cdot\,\,
\chi_{_{k\m(zz'0)}}(\tilde c_{z,z'})^{^{\sk}}\,
\chi_{_{(zz'0)\m(zz'z''0)}}(\tilde c_{z,z'})^{^{\sk}}\,
\chi_{zz'z''0}((\delta\tilde c)_{z,z',z''})^{^{\sk}}\,.
\qqq
The first factor may be rewritten as 
$\,\chi_{_{l\m(zz'z''0)}}(w_z(\tilde c_{z',z''})
w_z^{-1})^{^{-\sk}}\,$ using the 2$^{\rm nd}$ identity 
in (\ref{toc}) and combines with the next two to
\qq
&\displaystyle{\chi_{_{l\m(zz'z''z''0)}}((\delta\tilde c)_{z,z',z''})^{^{-\sk}}
\,\,\chi_{_{l\m(zz'z''z''0)}}(\tilde c_{z,z'})^{^{-\sk}}}&\cr\cr
&\displaystyle{=\ \chi_{_l}((\delta\tilde c)_{z,z',z''})^{^{\sk}}\,
\chi_{_{zz'z''0}}((\delta\tilde c)_{z,z',z''})^{^{-\sk}}\,
\chi_{_{l\m(zz'z''z''0)}}(\tilde c_{z,z'})^{^{-\sk}}\,.}&
\qqq
With the next three factors, it reproduces with the use of property
(\ref{rel3}) the right hand side of (\ref{rtc}). 
\smallskip

\hspace{10cm}$\Box$
\vskip 0.4cm

\noindent{\bf Proof of Lemma 2.}\ \ This proceeds similarly.
With the use of the explicit 
expression (\ref{rtc1}) the middle term of (\ref{tcc}) becomes
\qq
&&\,\,\,\chi_{_{(z'z''0)\m(z'z''z'''0)}}(\tilde c_{z',z''})^{^{\sk}}\,
\,\chi_{_{(zz'z''0)\m(zz'z''z'''0)}}(\tilde c_{zz',z''})^{^{-\sk}}\cr\cr
&&\cdot\,\,\chi_{_{(zz'z''0)\m(zz'z''z'''0)}}(\tilde c_{z,z'z''})^{^{\sk}}
\,\,\chi_{_{(zz'0)\m(zz'z''z'''0)}}(\tilde c_{z,z'})^{^{-\sk}}\cr\cr
&&\cdot\,\,\chi_{_{(zz'0)\m(zz'z''0)}}(\tilde c_{z,z'})^{^{\sk}}
\,\,\chi_{_{z'z''z'''0}}((\delta\tilde c)_{z',z'',z'''})^{^{\sk}}\cr\cr
&&\cdot\,\,\chi_{_{zz'z''z'''0}}((\delta\tilde c)_{zz',z'',z'''})^{^{-\sk}}
\,\,\chi_{_{zz'z''z'''0}}((\delta\tilde c)_{z,z'z'',z'''})^{^{\sk}}\cr\cr
&&\cdot\,\,\chi_{_{zz'z''z'''0}}((\delta\tilde c)_{z,z',z''z'''})^{^{-\sk}} 
\,\,\chi_{_{zz'z''0}}((\delta\tilde c)_{z,z',z''})^{^{\sk}}\,.
\qqq
The first factor is equal to  
$\,\chi_{_{(zz'z''0)\m(zz'z''z'''0)}}(w_z\,\tilde c_{z',z''}\m
w_z^{-1})^{^{\sk}}$, \,see (\ref{toc}), and it combines with the next 
four ones to
\qq
\chi_{_{(zz'z''0)\m(zz'z''z'''0)}}((\delta\tilde c)_{z,z',z''})^{^{\sk}}
\ =\ \chi_{_{zz'z''0}}((\delta\tilde c)_{z,z',z''})^{^{-\sk}}
\,\chi_{_{zz'z''z'''0}}((\delta\tilde c)_{z,z',z''})^{^{\sk}}\,,
\qqq
see (\ref{rel3}) and (\ref{rel2}). Together with the 
remaining factors, one obtains, rewriting the sixth factor 
as $\,\chi_{_{zz'z''z'''0}}(w_z(\delta\tilde c_{z',z'',z'''})
w_z^{-1})^{^{\sk}}$, \,an expression that reduces to 
\qq
\chi_{_{zz'z''z'''0}}((\delta^2c)_{z,z',z'',z'''})^{^{\sk}}
\qqq
and is equal to $\,1\,$ due to triviality of $\,\delta^2$.
\medskip

\hspace{10cm}$\Box$

\end{document}